          \documentclass{cmslatex}
\usepackage[paperwidth=7in, paperheight=10in, margin=.875in]{geometry}
          \usepackage{graphics}
          \usepackage{graphicx}
          \usepackage{amssymb}

\renewcommand{\theequation}{\arabic{section}.\arabic{equation}}

          \newtheorem{thm}{Theorem}[section]
          \newtheorem{prop}[thm]{Proposition}
          \newtheorem{lem}[thm]{Lemma}
          \newtheorem{cor}[thm]{Corollary}

          \newtheorem{rem}[thm]{Remark}

          \newcommand{\phie}{\phi_\epsilon}
          
          \newcommand{\sumk}{\sum_{k=1}^{N_1}}
          \newcommand{\suml}{\sum_{l=1}^{N_2}}
          \newcommand{\eapx}{e^{a_k\phie(x)}}
          \newcommand{\ebpx}{e^{-b_l\phie(x)}}
          \newcommand{\eapy}{e^{a_k\phie(y)}dy}
          \newcommand{\ebpy}{e^{-b_l\phie(y)}dy}
          \newcommand{\intt}{\int_{-1}^1}
          \newcommand{\etae}{\eta_\epsilon}
          \newcommand{\ds}{\displaystyle}
          \newcommand{\BE}{\begin{equation}}
          \newcommand{\EE}{\end{equation}}

\begin{document}
          \title{Boundary layer solutions of Charge Conserving Poisson-Boltzmann equations: one-dimensional case\thanks{The work is done when C.-C. Lee and T.-C. Lin
were visiting the Department of Mathematics of Penn State University (PSU) during the summer of 2014. They express sincere thanks to PSU for
all the hospitality and the great working environment. The authors also thank Professors Bob
Eisenberg, Ping Sheng, Chiun-Chuan Chen and Rolf
Ryham for numerous interesting discussions and valuable comments.
The authors also thank the anonymous referees for their helpful comments and suggestions.\newline
The research of C.-C. Lee was partially supported by the Ministry of Science and Technology of Taiwan under the grants NSC-101-2115-M-134-007-MY2 and MOST-103-2115-M-134-001. Y. Hyon is partially supported by the National Institute for Mathematical Sciences grant funded by the Korean government (No. B21501). T.-C. Lin is partially supported by National Center of Theoretical Sciences (NCTS) and the Ministry of Science and Technology of Taiwan grant NSC-102-2115-M-002-015-MY4 and MOST-103-2115-M-002-005-MY3. C. Liu is partially supported by the NSF grants DMS-1109107, DMS-1216938,  DMS-1159937 and DMS-1412005.}}

          \author{Chiun-Chang Lee\thanks{ Department of Applied Mathematics, National Hsinchu University of Education, No. 521, Nanda Road, Hsinchu 300, Taiwan,  email(chlee@mail.nhcue.edu.tw)}\and{Hijin Lee\thanks{ Department of Mathematical Sciences, Korea Advanced Institute of Science and Technology, Daejeon, Republic of Korea 305-701,  email(hijin@kaist.ac.kr) }}\and{YunKyong Hyon\thanks{ Division of Mathematical Models, National Institute for Mathematical Sciences, Daejeon, Republic of Korea 305-811,  email(hyon@nims.re.kr) }}\and{Tai-Chia Lin\thanks{ Department of Mathematics, Center for Advanced Study in Theoretical Sciences (CASTS), National Center for Theoretical Sciences (NCTS), National Taiwan University, Taipei, Taiwan 10617,  email(tclin@math.ntu.edu.tw) }}\and{Chun Liu\thanks{ Department of Mathematics, Pennsylvania State University, University Park,  PA 16802, USA,  email(liu@math.psu.edu) }}}





         \pagestyle{myheadings} \markboth{Boundary layer solutions of Charge Conserving Poisson-Boltzmann equations}{C.-C. Lee, H. Lee, Y. Hyon. T.C. Lin and C. Liu} \maketitle

          \begin{abstract}
          For multispecies ions, we study boundary layer solutions of charge conserving Poisson-Boltzmann (CCPB) equations~\cite{WXL13} (with a small parameter $\epsilon$) over a finite one-dimensional (1D) spatial domain, subjected to Robin type boundary conditions with variable coefficients. Hereafter, 1D boundary layer solutions mean that as $\epsilon$ approaches zero, the profiles of solutions form boundary layers near boundary points and become flat in the interior domain. These solutions are related to electric double layers with many applications in biology and physics. We rigorously prove the asymptotic behaviors of 1D boundary layer solutions at interior and boundary points. The asymptotic limits of the solution values (electric potentials) at interior and boundary points with a potential gap (related to zeta potential) are uniquely determined by explicit nonlinear formulas (cannot be found in classical Poisson-Boltzmann equations) which are solvable by numerical computations.
          \end{abstract}

          \section{Introduction}\label{intro}
Almost all biological activities involve transport in ionic solutions, which involves
various couplings and interactions of multiple species of ions. Many complicated types of electrolytes involved in biological processes, such as those in ion channel proteins, certain amino acids (movable side chain) are crucial to the functions of these ion channels. The electrostatic properties involving multispecies (at least three species) ions can be fundamentally different to those with only one or two species~\cite{BVHEG07, Lw}.
To see such difference, we study charge conserving Poisson-Boltzmann (CCPB) equation for multispecies ions which is derived from steady state Poisson-Nernst-Planck systems with charge conservation law, and is the surface potential model for the generation of a surface charge density layer related to electric double layers~\cite{LHLL, WXL13}. For simplicity of analysis, we consider a physical domain $x\in (-1,1)$ with the simplest geometry, and represent CCPB equation as follows:
\begin{equation}\label{p4} -\epsilon^2{\phi }''=\sum\limits_{i=1}^{N}{\frac{{{z}_{i}}{{e}_{0}}{{m}_{i}}}{\int_{-1}^{1}{{{e}^{-\frac{{{z}_{i}}{{e}_{0}}}{{{k}_{B}}T}\phi }}}}{{e}^{-\frac{{{z}_{i}}{{e}_{0}}}{{{k}_{B}}T}\phi }}}\quad \mbox{ for}\quad x\in \left( -1,1 \right)\,,
\end{equation}
where $m_i$ is the total concentration of species $i$ with the valence $z_i$, $\phi$ is the (electrical) potential, $e_0$ is the elementary charge, $k_B$ is the Boltzmann constant, and $T$ is the absolute temperature. The
parameter ~$\epsilon=\left(\epsilon_0\,U_T/(d^2\,e\,S)\right)^{1/2}>0$,
where $\epsilon_0$ is the dielectric constant of the electrolyte,
$U_T$ is the thermal voltage, $d$ is the length of the domain
$(-1,1)$, and $S$ is the appropriate concentration scale
(cf.~\cite{PJ}). Furthermore, $\epsilon\,d$ is known as the Debye
length and $\epsilon$ is of order $10^{-2}$ for the physiological
cases of interest (cf.~\cite{BCE}). Thus we may assume $\epsilon$
as a small parameter tending to zero. Similar equations to (\ref{p4}) can also be obtained by the other variational method~\cite{ZWL}.

Under suitable scales on $\phie$ and $\epsilon$, we let $-a_i$'s be the valences of anions, i.e., $a_k= -z_k$, $k=1,\cdots,N_1$ and $b_j$'s be the valences of cations, i.e., $b_l=z_l$, $l=1,\cdots,N_2$. Then the total concentrations of anions and cations are approximately given as $\alpha_k\sim m_k$ ($k=1,\cdots,N_1$) and $\beta_l\sim m_l$ ($l=1,\cdots,N_2$), respectively. Hence equation~(\ref{p4}) can be transformed into
\begin{eqnarray}
 \epsilon^2 \phie''(x) &=& \sumk \frac{a_k \alpha_k}{\intt \eapy}\,\eapx -\suml \frac{b_l \beta_l}{\intt \ebpy}\,\ebpx\nonumber \\[-0.5em]
 &&\label{2-eqn1}\\[-0.5em]
 &&\quad\hbox{for}\quad x \in (-1,1) \,,\nonumber
\end{eqnarray}
where $a_k$'s and $b_l$'s satisfy $1\leq a_1<a_2<\cdots<a_{N_1}$ and $1\leq b_1<b_2<\cdots<b_{N_2}$.

Most of the physical and biological systems possess the charge neutrality (zero net charge). One may assume the pointwise charge neutrality i.e.~at all points the anion and cation charges exactly cancel in order to make calculations easier in a free diffusion system (cf.~\cite{H1} p.~319). Here we replace the pointwise charge neutrality by a weaker hypothesis called the global electroneutrality being represented as
\begin{equation}\label{g-neutral}
{\mbox{Global Electro-neutrality: }\sum_{k=1}^{N_1}\,a_k\alpha_k=\sum_{l=1}^{N_2}\,b_l\beta_l\,,}
\end{equation}
which means that the total charges of anions and cations are equal, where $-a_k$'s and $b_l$'s are the valences, and $\alpha_k$'s and $\beta_l$'s are the concentrations of anions and cations, respectively. Consequently, the CCPB equation (\ref{2-eqn1}) may satisfy (\ref{g-neutral}).

When one deals with more general (realistic) situations, such as when there are more than two species involved in the solution, situations become more subtle and complicated. Note that the equation~(\ref{2-eqn1}) has nonlocal dependence on $\phie$ which is essentially different from the classical Poisson-Boltzmann (PB) equation as follows:
\begin{equation}\label{2-eqnw}
\epsilon^2\phie''(x)=\sum_{k=1}^{N_1}\frac{a_k\alpha_k}{2}e^{a_k\phie(x)}-\sum_{l=1}^{N_2}\frac{b_l\beta_l}{2}e^{-b_l\phie(x)}\quad\hbox{for}\quad x\in(-1,1)\,.
\end{equation}
Here $\frac{\alpha_k}{2}$'s and $\frac{\beta_l}{2}$'s are bulk concentration of anions and cations, respectively. In equation (\ref{2-eqn1}), $\alpha_k$'s and $\beta_l$'s are for total concentration of anions and cations, respectively. For notation convenience, we use the same notations $\alpha_k$'s and $\beta_l$'s in equations (\ref{2-eqn1}) and (\ref{2-eqnw}), but with different physical meaning. In this paper, we shall show different asymptotic behaviors of the CCPB equation~(\ref{2-eqn1}) and the PB equation~(\ref{2-eqnw}) for various constants $N_1, N_2, a_k, \alpha_k, b_l, \beta_l$ satisfying (\ref{g-neutral}). The main goal of this paper is to compare the CCPB equation~(\ref{2-eqn1}) and the PB equation~(\ref{2-eqnw}) under the hypothesis (\ref{g-neutral}). Such a difference can be clarified in Theorems~\ref{2-thm:w} and \ref{2-thm5}, see also, Remark~\ref{pbn-pb}.

Boundary effects are important in a wide range of applications and provide formidable challenges \cite{Hunter1,I}. For CCPB equations, the main issue is how boundary conditions effect the solution values (electric potentials) at interior and boundary points. One may use the Neumann boundary condition for a given surface charge distribution and the Dirichlet boundary condition for a given surface potential (cf.~\cite{A1}). Here we consider a Robin boundary condition~\cite{BCB, HyonMori, LHLL, LMB, R1,RLZ,RLW, ZMB} for the electrostatic potential $\phi$ at $x=\pm 1$ is given by
\begin{equation}\label{2-eqn2}
\phie(1)  +\etae \phie'(1)  = \phi_0^+, \quad \phie(-1) -\etae \phie'(-1) = \phi_0^-\,.
\end{equation}
where $\phi_0^+$, $\phi_0^-$ are extrachannel electrostatic potentials and $\eta_\epsilon\geq 0$ is the coefficient depending on the dielectric constant \cite{Mori1, Mori2}, and related to the surface capacitance. The parameter ratio $\eta_\epsilon=\epsilon_S/C_S$ can be viewed as a measure of the Stern layer thickness, where $\epsilon_S$ and $C_S$ are the effective permittivity and the capacitance of the Stern layer, respectively (cf. \cite{BCB}). Thus we may regard $\frac{\eta_\epsilon}{\epsilon}$ as the ratio of the Stern-layer width to the Debye screening length. Similar discussion can also be found in~\cite{DC} and~\cite{O}. To see the influence of $\frac{\eta_\epsilon}{\epsilon}$ on the asymptotic behavior of $\phie$'s, we consider the limit $\lim_{\epsilon\downarrow0}\frac{\eta_\epsilon}{\epsilon}$ to be either a non-negative constant $\gamma$ or infinity.

\begin{figure}[hbt]
 \centering{%
  \begin{tabular}{@{\hspace{-0pc}}c}
  \includegraphics[width=3.8in]{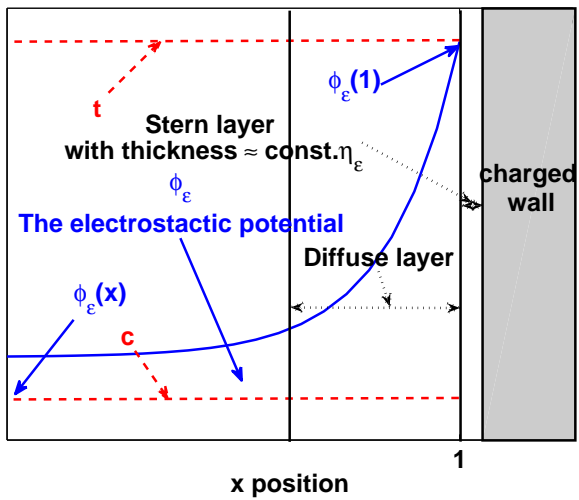}
  \end{tabular}}
 \caption{\small Schematic picture of Robin boundary condition, $\phi_{\epsilon} \pm \eta_{\epsilon}(\phi_{\epsilon})_x = \phi_0^{\pm}$ at $x=\pm 1$, and the limit values $t=\lim_{\epsilon\to 0}\phi_\epsilon(1)$, $c=\lim_{\epsilon\to 0}\phi_\epsilon(x)$, $x\in (-1,1)$ and $\zeta=t-c$.}\label{fig:layer}
\end{figure}

Suppose $\lim_{\epsilon\downarrow0}\frac{\eta_\epsilon}{\epsilon}=\gamma$ i.e.~$\eta_\epsilon\sim \gamma\epsilon$, where $\gamma$ is a non-negative constant. Then we show that the solution $\phie$ of (\ref{2-eqn1}) with (\ref{2-eqn2}) satisfies $\lim_{\epsilon\downarrow0}\phie(\pm 1)=\pm t$ and $\lim_{\epsilon\downarrow0}\phie(x)=c$ for $x\in(-1,1)$, where $c$ and $t$ can be uniquely determined by (\ref{1.19.1})-(\ref{1.19.3}) which imply that the value $c$ is changed with respect to $t$ (see Figure~\ref{fig:layer}). Moreover, the potential difference $\zeta=t-c$ is decreasing to $\gamma$ (cf.~Theorem~\ref{2-thm5}). Note that as the parameter $\epsilon$ goes to zero, the solution~$\phi_\epsilon$ has a boundary layer producing the potential gap $\zeta=t-c$ affected by Stern and Debye (diffuse) layers and related to zeta potential (cf.~\cite{Hunter}) which plays an important role in ionic fluids. However, for the PB equation~(\ref{2-eqnw}), the value $c$ must be zero which is independent of $t$ and $\gamma$ (cf.~Theorem~\ref{2-thm:w}). This shows the difference of the CCPB equation~(\ref{2-eqn1}) and the PB equation~(\ref{2-eqnw}) which can also be observed by numerical experiments (See Figure~\ref{fig:F_1} and Table \ref{table:D_c} in Section~\ref{2-mpnpsec:5}). Furthermore, numerical computations give several conditions to let the profile of function $c$ to $\gamma$ become monotone decreasing and increasing (Figure~\ref{fig:ct_3_12_allm} and~\ref{fig:ct_4_all} in Section~\ref{2-mpnpsec:5}) and non-monotone (Figure~\ref{fig:ct_5_14_all_z} in Section~\ref{2-mpnpsec:5}).

In \cite{LHLL}, we studied the CCPB equation~(\ref{2-eqn1}) for case of $N_1=N_2=1$, $(a_1,b_1)=(1,1)$ and $(\alpha_1,\beta_1)=(\alpha,\beta)$, i.e.,
 the case of one  anion and one cation species with monovalence.
In this case, equation (\ref{2-eqn1}) can be rewritten as
\begin{eqnarray}
 \epsilon^2 \phie''(x) = n_\epsilon(x) -p_\epsilon(x)&&
 \quad\hbox{for}\quad x \in (-1,1),\label{2-eqn2014}\\
 n_\epsilon(x)\equiv\frac{\alpha e^{\phie(x)}}{\int_{-1}^1e^{\phie(y)} dy}\quad\mathrm{and}&&\quad  p_\epsilon(x)\equiv\frac{\beta e^{-\phie(x)}}{\int_{-1}^1e^{-\phie(y)} dy}, \label{2013-1013-3}
\end{eqnarray}
where $ n_\epsilon(x)$ and $p_\epsilon(x)$
represent (pointwise) concentrations of anion and cation species, respectively.
When $\alpha=\beta$ holds (the electroneutral case), we had shown previously
that $\lim_{\epsilon\downarrow0}n_\epsilon(x)=\lim_{\epsilon\downarrow0}p_\epsilon(x)=\frac{\alpha}{2}$ for $x\in (-1,1)$. Moreover, the CCPB equation~(\ref{2-eqn2014})-(\ref{2013-1013-3}) and the conventional PB equation
$\epsilon^2w_\epsilon''(x)=\frac{\alpha}{2}\left( e^{w_\epsilon(x)}- e^{-w_\epsilon(x)}\right)$ have same asymptotic behavior (cf. Theorem 1.4 of \cite{LHLL}). In order for the readers to compare those with the results in the current paper, most results of~\cite{LHLL} are summarized in Appendix. To certain degrees, it also justifies why in many situations, PB equation provides more or less expected solutions. On the other hand, we consider the non-electroneutral case, i.e. $\alpha\neq\beta$. Without loss of generality, we assume $\alpha<\beta$ i.e. $\int_{-1}^1n_\epsilon(x)dx<\int_{-1}^1p_\epsilon(x)dx$ which means that the total concentration of anion species is less than that of cation species. Then we prove that
$\lim_{\epsilon\downarrow0}n_\epsilon(x)=\lim_{\epsilon\downarrow0}p_\epsilon(x)=\frac{\alpha}{2}$ for $x\in(-1,1)$, but $\lim_{\epsilon\downarrow0}n_\epsilon(\pm1)=0<\lim_{\epsilon\downarrow0}\epsilon^2p_\epsilon(\pm1)=\frac{(\alpha-\beta)^2}{8}$
(cf. (\ref{2014-0910-2})). This shows that electroneutrality holds true in the interior of $(-1,1)$, but non-electroneutrality occurs at the boundary points $\pm 1$. Furthermore, the extra charges are accumulated near the boundary points $\pm 1$ (see Theorem~~\ref{NE-thm}).

The mixture of monovalent and divalent ions such as $\rm{Na^+}$, $\rm{K^+}$, $\rm{Cl^-}$ and $\rm{Ca^{2+}}$ plays the most important roles for vital biological processes. For instance, opening and closing of ionic channels is accomplished by escape or entry of $\rm{Ca^{2+}}$ into the channels (cf.~\cite{GW}). The voltage may depend on $[\rm{Ca^{2+}}]$ the concentration of $\rm{Ca^{2+}}$ (cf.~\cite{H1}). Differences in ionic concentrations create a potential gap across the cell membrane that drives ionic currents (cf.~\cite{K0}~P.~34). To see how the voltage i.e.~(electrical) potential depends on $[\rm{Ca^{2+}}]$, we may use the equation~(\ref{2-eqn1}) with $N_1=1$, $N_2=2$, $a_1=b_1=1$ and $b_2=2$ to describe the mixture of $\rm{Na^+}$ (or $\rm{K^+}$), $\rm{Cl^-}$ and $\rm{Ca^{2+}}$ ions, where $\alpha_1\sim [\rm{Cl^-}]$, $\beta_1\sim [\rm{Na^+}]$ and $\beta_2\sim [\rm{Ca^{2+}}]$. In Theorem~\ref{2-thm5}~(ii), we prove that when the electro-neutrality holds, that is, $\alpha_1=\beta_1+2\beta_2$, the solution $\phie$ of~(\ref{2-eqn1}) satisfies $\ds\lim_{\epsilon\to 0}\phie(x)=c$ for $x\in (-1,1)$ and $c\in (c_*, 0)$ is uniquely determined by (\ref{1.19.1}) and
\begin{equation}\label{ca1}
\frac{1-e^{3c}\cosh t}{e^c\,\sinh
c}=\frac{\beta_1}{\beta_2}=\frac{[\rm{Na^+}]}{[\rm{Ca^{2+}}]}\,>0\,
\end{equation}
where $t=\ds\lim_{\epsilon\to 0}\phie(1)>0$, and $c_*=\frac{1}{3}\log {\rm sech}\, t$ is a negative constant (see Remark~\ref{1.4}). The formula~(\ref{ca1}) shows that the interior potential (voltage) $c$ is increased if the boundary potential $t$ is fixed and the ratio $[\rm{Na^+}]/[\rm{Ca^{2+}}]$ is increased e.g.~$[\rm{Ca^{2+}}]$ is decreased and $[\rm{Na^+}]$ is fixed. Furthermore, Theorem~\ref{2-thm5} is also applicable to the other cases with multi-species ions including multivalent and polyvalent ions so the formula~(\ref{ca1}) can be generalized to
\begin{equation}\label{ca2} \frac{z -e^{(
1 +z ) c}( \frac{\sinh(zt)}{\sinh t})}{2 e^c \sinh
c}=\frac{\beta_1}{\beta_2}\,,
\end{equation}
for $a_1=b_1=1$ and $b_2=z\geq 2$ (see Remark~\ref{1.4}). Note that (\ref{ca2}) shows how the value $c$ depends on the value $t$. Such a result cannot be found in the PB equation~(\ref{2-eqnw}).

\subsection{Asymptotic behavior of the PB equation~(\ref{2-eqnw})-(\ref{2-eqn2})}
The PB equation (\ref{2-eqnw}) with the boundary condition
(\ref{2-eqn2}) can be regarded as the Euler-Lagrange equation of the energy functional
\begin{equation}\label{20140818-1}
E^{\mathrm{PB}}_\epsilon[u]=\frac{1}{2}\int_{-1}^1\left(\epsilon^2|u'|^2+f(u)\right)dx+
\frac{\epsilon^2}{2\eta_\epsilon}\left[(\phi_0^--u(-1))^2+(\phi_0^+-u(1))^2\right],
\end{equation}
for $u\in H^1((-1,1))$, where
\begin{equation}\label{id:2-015}
f(s)=\ds\sumk\alpha_ke^{a_ks}+\ds\suml\beta_le^{-b_ls}\quad\hbox{
for}\quad s\in\mathbb{R}.
\end{equation}

For the PB equation (\ref{2-eqnw}) with the boundary condition
(\ref{2-eqn2}), we study the asymptotic behavior of the
solution $\phie$ of (\ref{2-eqnw}) as $\epsilon$ approaches zero.
The boundary condition (\ref{2-eqn2}) plays a crucial role on
the monotonicity of $\phie$. Here we consider three cases for the
signs of $\phi_0^+$ and $\phi_0^-$:
(a)~$\min\{\phi_0^+,\phi_0^-\}>0$,
(b)~$\max\{\phi_0^+,\phi_0^-\}<0$ and
(c)~$\min\{\phi_0^+,\phi_0^-\}\leq0\leq\max\{\phi_0^+,\phi_0^-\}$.
Then the corresponding results are stated as follows:\\
\begin{theorem}\label{2-thm:w}
Assume $\sumk\,a_k\alpha_k=\suml\,b_l\beta_l$. Let $\phie\in
C^\infty((-1,1))\cap C^2([-1,1])$ be the solution of equation
(\ref{2-eqnw}) with the boundary condition (\ref{2-eqn2}). Then
\begin{enumerate}
\item[(i)] For $x\in(-1,1)$, $|\phie(x)|$ exponentially converges
to zero as $\epsilon$ goes to zero; \item[(ii)] If
$\min\{\phi_0^+,\phi_0^-\}>0$, then $\phie$ is convex on
$[-1,1]$ and $0\leq\phie(x)\leq\max\{\phi_0^+,\phi_0^-\}$ for $x\in[-1,1]$.
Moreover, there exists $\epsilon^*>0$ such that for
$0<\epsilon<\epsilon^*$,
 $\phie$ attains the minimum at an interior point of $(-1,1)$.
\item[(iii)] If $\max\{\phi_0^+,\phi_0^-\}<0$, then $\phie$ is
concave on $[-1,1]$ and
$\min\{\phi_0^+,\phi_0^-\}\leq\phie(x)\leq0$ for $x\in[-1,1]$.
Moreover, there exists $\epsilon^*>0$ such that for
$0<\epsilon<\epsilon^*$, $\phie$ attains the maximum at an
interior point of $(-1,1)$. \item[(iv)] If
$\min\{\phi_0^+,\phi_0^-\}\leq0\leq\max\{\phi_0^+,\phi_0^-\}$,
then $\phie$ is monotone on $[-1,1]$ and
$\min\{\phi_0^+,\phi_0^-\}\leq\phie\leq\max\{\phi_0^+,\phi_0^-\}$.
\item[(v)] If
$\displaystyle\lim_{\epsilon\downarrow0}\frac{\eta_\epsilon}{\epsilon}=\gamma$
and $0\leq\gamma<\infty$, then
$\displaystyle\lim_{\epsilon\downarrow0}\phie(1)=\widehat{t}$
uniquely determined by
\begin{equation}\label{1.25}
|\phi_0^+-\widehat{t}|=\gamma(f(\widehat{t})-f(0))^{1/2}\quad
and\quad \min\{0,\phi_0^+\}\leq \widehat{t} \leq
\max\{0,\phi_0^+\},
\end{equation}
where $f$ is defined by (\ref{id:2-015}).
Moreover, $\widehat{t}\equiv\widehat{t}(\gamma)$ is decreasing in
$\gamma$ if $\phi_0^+>0$ and increasing in $\gamma$ if
$\phi_0^+<0$.
\end{enumerate}
\end{theorem}

\subsection{The main results}
\medskip
\noindent
In this section we present the main results, which are about the asymptotic behavior
of the solution $\phi_\epsilon$ of (\ref{2-eqn1}) and (\ref{2-eqn2}) as
$\epsilon$ goes to zero, in our research of CCPB equation.
The  CCPB equation~(\ref{2-eqn1}) with the boundary condition~(\ref{2-eqn2})
can be regarded as the Euler-Lagrange equation of the energy functional
\begin{eqnarray}
E_\epsilon[u]&=&\int_{-1}^1\frac{\epsilon^2}{2}|u'|^2dx+\sumk\alpha_k\log\int_{-1}^1e^{a_ku}dx+\suml\beta_l\log\int_{-1}^1e^{-b_lu}dx\nonumber\\[-0.7em]
&&\label{2014-0821-energy}\\[-0.7em]
&&~~~~~~+\frac{\epsilon^2}{2\eta_\epsilon}\left[(\phi_0^--u(-1))^2+(\phi_0^+-u(1))^2\right],\nonumber
\end{eqnarray}
for $u\in H^1((-1,1))$. The existence and uniqueness for the solution of (\ref{2-eqn1})
and (\ref{2-eqn2}) is the following proposition:\\
\begin{proposition}\label{2-thm1}
 There exists a unique solution $\phie\in\,C^\infty((-1,1))\cap\,C^2([-1,1])$ of
 the equation (\ref{2-eqn1}) with the boundary condition (\ref{2-eqn2}).\\
\end{proposition}

\noindent The proof of the above Proposition~\ref{2-thm1}
 can be easily obtained from the arguments of \cite{LHLL} (see Appendix therein) and \cite{Lcc}.

Suppose $\phi_0^+=\phi_0^-=A$ and
$\sum_{k=1}^{N_1}\,a_k\alpha_k=\sum_{l=1}^{N_2}\,b_l\beta_l$. Then
Proposition~\ref{2-thm1} implies the solution of (\ref{2-eqn1})
and (\ref{2-eqn2}) must be trivial and $\phie\equiv A$. To study
the nontrivial solution of (\ref{2-eqn1}) and (\ref{2-eqn2}), it
is sufficient to assume $\phi_0^+\neq\phi_0^-$. Replacing
$\phie$ by $\phie+C$ for any constant $C$, one may remark that
the equation~(\ref{2-eqn1}) is invariant. Consequently, without
loss of generality, we may assume $-\phi_0^-=\phi_0^+>0$
hereafter.

When $\sum_{k=1}^{N_1}\,a_k\alpha_k=\sum_{l=1}^{N_2}\,b_l\beta_l$, i.e., the global electroneutral case,
Theorem~\ref{2-thm2} shows that $\max_{x\in[-1,1]}|\phie(x)|$ is uniformly bounded to $\epsilon$
and that  $\phie'$ exponentially approaches zero in $(-1,1)$ as $\epsilon$ tends to zero.
Thus, it is expected that there exists a constant $c$ such that all interior values of $\phie$ tends to $c$ as $\epsilon$ goes to zero.
Along with Lebesgue's dominated convergence theorem, we have
\begin{equation}\label{2014-0902}
\lim_{\epsilon\downarrow0}\int_{-1}^1e^{a_k\phie}dx=2e^{a_kc},\quad\lim_{\epsilon\downarrow0}\int_{-1}^1e^{-b_l\phie}dx=2^{-b_lc},
\end{equation}
and then the energy functional (\ref{2014-0821-energy}) with $u=\phie$ approaches to the energy functional $\widehat{E}_\epsilon[\phie]$ as follows (up to a constant independent of $\phie$):
\begin{equation}\label{2014-0821-PB}
~~~~~~\widehat{E}_\epsilon[\phie]=\frac{1}{2}\int_{-1}^1(\epsilon^2|\phie'|^2+f(\phie-c))dx+
\frac{\epsilon^2}{2\eta_\epsilon}\left[(\phi_0^--\phie(-1))^2+(\phi_0^+-\phie(1))^2\right],
\end{equation}
where $f$ is defined by (\ref{id:2-015}).
Here we have used $\lim_{\epsilon\downarrow0}\left(\frac{1}{2}\int_{-1}^1e^{a_k(\phie-c)}dx-1\right)=0$ (by (\ref{2014-0902})) and
the approximation $\log(1+s)\sim s$ with $s=\frac{1}{2}\int_{-1}^1e^{a_k(\phie-c)}dx-1$ to get
$$\log\int_{-1}^1e^{a_k\phie}dx\sim\frac{1}{2}\int_{-1}^1e^{a_k(\phie-c)}dx+\log(2e^{a_kc})-1\quad\mathrm{as}\quad  0<\epsilon\ll1.$$
Similarly, we have
$$\log\int_{-1}^1e^{-b_l\phie}dx\sim\frac{1}{2}\int_{-1}^1e^{-b_l(\phie-c)}dx+\log(2e^{-b_lc})-1\quad\mathrm{as}\quad 0<\epsilon\ll1.$$
Therefore, we show that in the case of global electroneutrality (\ref{g-neutral}), the energy functional~(\ref{2014-0821-energy}) approaches (\ref{2014-0821-PB}), which has the same form as the PB energy functional~(\ref{20140818-1}).

The asymptotic behavior of $\phie$'s at boundary $x=\pm1$ may depend on the scale of
$\eta_\epsilon$. Here we study two cases for the scale of
$\eta_\epsilon\geq 0$:
(i)
$\lim_{\epsilon\downarrow0}\frac{\eta_\epsilon}{\epsilon}=\infty$ and
(ii)
$\lim_{\epsilon\downarrow0}\frac{\eta_\epsilon}{\epsilon}=\gamma$,
where $\gamma$ is a nonnegative constant.
Then the relation between the boundary value limits
$\ds\lim_{\epsilon\to 0}\phie(\pm1)$ and the interior value limit
$c$ are demonstrated as
follows:\\
\begin{theorem}\label{2-thm5}
Assume $-\phi_0^-=\phi_0^+>0$ and
$\ds\sumk\,a_k\alpha_k=\ds\suml\,b_l\beta_l$. Let $\phie\in
C^\infty((-1,1))\cap C^2([-1,1])$ be the solution of equation
(\ref{2-eqn1}) with the boundary condition (\ref{2-eqn2}).Then
$$\lim_{\epsilon\downarrow0}\phie(-1)=-t,\quad\lim_{\epsilon\downarrow0}\phie(1)=t\quad
\hbox{ and }\quad\lim_{\epsilon\downarrow0}\phie(x)=c\quad\hbox{
for }\quad x\in(-1,1)\,, $$ where $t$ and $c$ are determined as
follows:
\begin{enumerate}
\item[(i)]~~~If~~
$\ds\lim_{\epsilon\downarrow0}\frac{\eta_\epsilon}{\epsilon}=\infty$,
then $c=t=0$.
\item[(ii)]~If~~$\ds\lim_{\epsilon\downarrow0}\frac{\eta_\epsilon}{\epsilon}=\gamma$
and $0\leq\gamma<\infty$, then $(t,c)$ uniquely solves the
following equations:
\begin{eqnarray}
\phi_0^+-t&=&\gamma(f(t-c)-f(0))^{1/2}, \label{1.19.1}\\
f(t-c)&=&f(-t-c), \label{1.19.2}\\
|c|<&t&\leq\,\phi_0^+. \label{1.19.3}
\end{eqnarray}
Moreover, writing $t=t(\gamma)$ and $c=c(\gamma)$ in (ii), we have
\begin{itemize}
\item[(A)] $\displaystyle\lim_{\gamma\rightarrow0}t(\gamma)=\phi_0^+$,
 $\displaystyle\lim_{\gamma\rightarrow0}c(\gamma)=c^*$ and $\displaystyle\lim_{\gamma\rightarrow\infty}t(\gamma)=\lim_{\gamma\rightarrow\infty}c(\gamma)=0$,
 where $|c^*|<\phi_0^+$ is uniquely determined by $f(\phi_0^+-c^*)=f(-\phi_0^+-c^*)$.
\item[(B)] $t(\gamma)$ and $t(\gamma)-c(\gamma)$ both are
decreasing on $(0,\infty)$.\\
\end{itemize}
\end{enumerate}
\end{theorem}

Formally, using ${{\phi }_{\epsilon }}\to c$ in $\left( -1,1 \right)$ as $\epsilon$ tends to zero, equation (\ref{2-eqn1}) may approach to the following PB equation:
\begin{equation}\label{2-eqn1-asp}
\epsilon^2\phie''(x)=\sum_{k=1}^{N_1}\frac{a_k\alpha_k}{2}e^{a_k\left(\phie(x)-c\right)}-\sum_{l=1}^{N_2}\frac{b_l\beta_l}{2}e^{-b_l\left(\phie(x)-c\right)}\quad\hbox{for}\quad x\in(-1,1)\,,
\end{equation}
which may give results of Theorem~\ref{2-thm5} by formal asymptotic analysis. However, in this paper, we focus on rigorous mathematical analysis and provide the proof of Theorem~\ref{2-thm5} in Section~\ref{2-mpnpsec:2}.

Theorem~\ref{2-thm5}(i) shows that
there is no boundary layer and $\phie\rightarrow 0$ uniformly in
$[-1,1]$ as $\epsilon\downarrow0$ if $\lim_{\epsilon\downarrow0}\frac{\eta_\epsilon}{\epsilon}=\infty$.
Theorem~\ref{2-thm5}(ii) assures the existence of
boundary layers. Furthermore, Theorem~\ref{2-thm5}~(ii-A) and (ii-B) represent
the ratio of Stern screening length to the Debye screening length affects the boundary and interior potentials:
(a) The decrease of $\gamma$ results in the increase of $t-c$ (the potential difference between the boundary and interior); (b) If $\gamma\rightarrow\infty$, the potential difference $t-c$ may approach zero.
Notice that the formula~(\ref{1.25}) is quite different from
(\ref{1.19.1})-(\ref{1.19.3}). This may show the difference
between solutions of the CCPB equation~(\ref{2-eqn1}) and the PB
equation~(\ref{2-eqnw}).\\

\begin{rem}\label{pbn-pb}
\begin{itemize}
\item[(a)] Theorem~\ref{2-thm:w}~(ii) and~(iii) show that as
$\phi_0^+\phi_0^->0$, the solution $\phie$ of the PB equation
(\ref{2-eqnw}) may lose the monotonicity. However, the solution of
the CCPB equation (\ref{2-eqn1}) always keeps the monotonicity
(see Remark~\ref{2-rk1.2}~(i)). This provides the difference
between solutions of the CCPB equation (\ref{2-eqn1}) and the PB
equation (\ref{2-eqnw}).
\item[(b)] For equation (\ref{2-eqn1}), the values $c$ (interior potential) and $t$ (boundary potential) depend on each other and satisfy precise formulas (\ref{1.19.1})-(\ref{1.19.3}). However, for equation (\ref{2-eqnw}), interior potential and boundary potentials (determined by (\ref{1.25})) are independent to each other.\\
\end{itemize}
\end{rem}

\begin{rem}\label{1.4}
When $N_1=1$, $N_2=2$, $a_1=b_1=1$, $b_2=2$ and
$\alpha_1=\beta_1+2\beta_2$, we may get~(\ref{ca1}) from
(\ref{1.19.1}) and (\ref{1.19.2}). Moreover, (\ref{ca2}) can also
be derived from (\ref{1.19.1}) and (\ref{1.19.2}) for the case
that $N_1=1$, $N_2=2$, $a_1=b_1=1$, $b_2=z\geq 2$ and
$\alpha_1=\beta_1+z\beta_2$. By~(\ref{ca1}) and (\ref{ca2}), it is
easy to check that $\frac{d c}{dt}<0$ for $t>0$. Then
 $c=c(t)$ can be regarded as an decreasing function to $t>0$. Consequently, by
Theorem~\ref{2-thm5}~(iv), $c$ is increasing to $\gamma$.\\
\end{rem}

When $N_1=1$, $N_2=2$, $a_1=b_1=1$
and $b_2=2$, further asymptotic behavior of $\phie$ near the
boundary $x=\pm1$ describing the boundary layers is stated as
follows:\\
\begin{theorem}\label{2-thm:v}
Assume $N_1=1$, $N_2=2$, $a_1=b_1=1$ and $b_2=2$. Under the same
hypotheses of Theorem~\ref{2-thm2} and Theorem~\ref{2-thm5}(ii),
 the asymptotic behavior of $\phie$ near the boundary
$x=\pm1$ can be represented by
\begin{eqnarray}\label{2-id:thmw1}
\phi_{1,\epsilon}^+(x)\leq\,\phie(x)\leq\,\phi_{2,\epsilon}^+(x)\quad &&for\quad x\in(y_\epsilon^+,1)\,, \\
\phi_{1,\epsilon}^-(x)\leq\,\phie(x)\leq\,\phi_{2,\epsilon}^-(x)\quad
&&for\quad x\in(-1,y_\epsilon^-)\label{2-id:thmv}\,,
\end{eqnarray}
where $-1<y_\epsilon^-<y_\epsilon^+<1$ satisfy
$\ds\lim_{\epsilon\downarrow0}v(y_\epsilon^\pm)=0$, and
\begin{eqnarray}
\phi_{i,\epsilon}^+(x)&=&
c+\log\left\{A_{i,\epsilon}^++B_{i,\epsilon}^+
\mathrm{csch}^{2}\left[\frac{C_{i,\epsilon}^+}{\epsilon}(1-x)+\log\,D_{i,\epsilon}^+\right]\right\},\label{2-id:w13}\\
\phi_{i,\epsilon}^-(x)&=&
c+\log\left\{A_{i,\epsilon}^--B_{i,\epsilon}^-
\mathrm{sech}^{2}\left[\frac{C_{i,\epsilon}^-}{\epsilon}(1+x)+\log\,D_{i,\epsilon}^-\right]\right\},\label{2-id:w13add}\quad\,
i=1,2\,.
\end{eqnarray}
Here $A_{i,\epsilon}^\pm$, $B_{i,\epsilon}^\pm$,
$C_{i,\epsilon}^\pm$ and $D_{i,\epsilon}^\pm$, $i=1,2$, are
constants depending on $\epsilon$ such that
$A_{i,\epsilon}^\pm\rightarrow1$,
$B_{i,\epsilon}^\pm\rightarrow1+\frac{\beta_2}{\alpha_1}$,
$C_{i,\epsilon}^\pm\rightarrow\sqrt{\alpha_1+\beta_2}$ and
$D_{i,\epsilon}^\pm\rightarrow\frac{\sqrt{\alpha_1e^{\pm\,t-c}+\beta_2}+\sqrt{\alpha_1+\beta_2}}{\sqrt{\alpha_1e^{\pm\,t-c}+\beta_2}-\sqrt{\alpha_1+\beta_2}}$
as $\epsilon$ goes to zero.\\
\end{theorem}

\noindent In the case of $N_1=1$, $N_2=2$, $a_1=b_1=1$ and $b_2=2$, we may solve equation (\ref{2-eqn1-asp}) precisely and get the form of (\ref{2-id:w13}) and (\ref{2-id:w13add}) near $x=1$ and $x=-1$, respectively. One may remark how the values $c, t, \alpha_1$ and $\beta_2$ affect the asymptotic behavior of $\phie$ near the boundary $x=\pm1$.

When $\alpha\neq\beta$ (the non-electroneutral case), the asymptotic behavior for the solution $\phie$, $n_\epsilon$ and $p_\epsilon$ of the equation~(\ref{2-eqn2014})-(\ref{2013-1013-3})
with the boundary condition~(\ref{2-eqn2}) is stated as follows:\\
\begin{theorem}\label{NE-thm}
Assume $0<\alpha<\beta$ and $\phi_0^-=\phi_0^+$. Let $\phie\in
C^\infty((-1,1))\cap C^2([-1,1])$ be the solution of the equation~(\ref{2-eqn2014})-(\ref{2013-1013-3}) with the boundary condition~(\ref{2-eqn2}) and $\eta_\epsilon\geq 0$. Then
\begin{itemize}
\item[(i)] When $0<\epsilon<1$ and $0<\kappa<1$,
there exists a positive constant $\lambda_\epsilon(\kappa)$ depending on $\epsilon$ and $\kappa$ such that $\lim_{\epsilon\downarrow0}\lambda_\epsilon(\kappa)=0$ and
\begin{equation}\label{2014-0810-1pm}
\frac{\alpha}{2}-\lambda_\epsilon(\kappa)\leq n_\epsilon(x)\leq p_\epsilon(x)\leq \frac{\beta}{2}+\lambda_\epsilon(\kappa),\quad\mathrm{for}\quad x\in[-1+\epsilon^\kappa,1-\epsilon^\kappa].
\end{equation}
Moreover, we have
\begin{eqnarray}
\lim_{\epsilon\downarrow0}n_\epsilon(\pm1)=0\quad \mathrm{and}\quad&&\lim_{\epsilon\downarrow0}\epsilon^2p_\epsilon(\pm1)
=\frac{(\alpha-\beta)^2}{8},\label{2014-0910-2}\\
\lim_{\epsilon\downarrow0}\sup_{x\in[-1+\epsilon^\kappa,1-\epsilon^\kappa]}\left|n_\epsilon(x)-\frac{\alpha}{2}\right|
&=&\lim_{\epsilon\downarrow0}\sup_{x\in[-1+\epsilon^\kappa,1-\epsilon^\kappa]}\left|p_\epsilon(x)-\frac{\alpha}{2}\right|=0,\label{thm2014-0829}\\
\lim_{\epsilon\downarrow0}\int_{-1}^{-1+\epsilon^\kappa}n_\epsilon(x)dx&=&\lim_{\epsilon\downarrow0}\int^{1}_{1-\epsilon^\kappa}n_\epsilon(x)dx=0,\label{thm1-6}\\
\lim_{\epsilon\downarrow0}\int_{-1}^{-1+\epsilon^\kappa}p_\epsilon(x)dx&=&\lim_{\epsilon\downarrow0}\int^{1}_{1-\epsilon^\kappa}p_\epsilon(x)dx=\frac{\beta-\alpha}{2}.\label{thm1-7}
\end{eqnarray}
\item[(ii)] Let $K$ be any compact subset of $(-1,1)$. When $0<\epsilon\ll1$ is sufficiently small,
the asymptotic expansion of $\phie(x)-\phie(\pm1)$ in $\epsilon$ with the exact leading-order term $\log\frac{1}{\epsilon^2}$ and second-order term $O(1)$
is given as follows:
\begin{equation}\label{2013-1015-1}
\phie(x)-\phie(\pm1)=\log\frac{1}{\epsilon^2}+\log\frac{(\alpha-\beta)^2}{4\alpha}+o_\epsilon(1),\quad\mathrm{for}\quad x\in K,
\end{equation}
where $o_\epsilon(1)$ denotes as a small quantity tending to zero as $\epsilon$ goes to zero.\\
\end{itemize}
Similar results also hold for $0<\beta<\alpha$.\\
\end{theorem}

\begin{rem}\label{newrk1-8}
\begin{itemize}
\item[(i)] To exclude the boundary layer of $\phie$ with thickness $\epsilon^2$ (cf. Theorem~1.6 of \cite{LHLL}),
we consider integrals of $n_\epsilon$ and $p_\epsilon$ over the interval $[-1+\epsilon^\kappa,1-\epsilon^\kappa]$,
where $0<\kappa<1$ is independent of $\epsilon$. Note that $n_\epsilon$ and $p_\epsilon$ can be represented by $\phie$ (see (\ref{2013-1013-3})),
and that Theorem~\ref{NE-thm}(ii) implies
$\lim_{\epsilon\downarrow0}\int_{-1+\epsilon^\kappa}^{1-\epsilon^\kappa}n_\epsilon(x)dx=\int_{-1+\epsilon^\kappa}^{1-\epsilon^\kappa}p_\epsilon(x)dx=\alpha>0$,
$\lim_{\epsilon\downarrow0}\left(\int_{-1}^{-1+\epsilon^\kappa}+\int^{1}_{1-\epsilon^\kappa}\right)n_\epsilon(x)dx=0$
and $\lim_{\epsilon\downarrow0}\left(\int_{-1}^{-1+\epsilon^\kappa}+\int^{1}_{1-\epsilon^\kappa}\right)p_\epsilon(x)dx=\beta-\alpha>0$.
This shows that as $\epsilon$ approaches zero,
both the total concentrations of anion and cation species in the bulk $[-1+\epsilon^\kappa,1-\epsilon^\kappa]$ tend to the same positive constant $\alpha$,
while the total concentrations of anion and cation species in the region $[-1,-1+\epsilon^\kappa)\cup(1-\epsilon^\kappa,1]$ ( which is next to the boundary
with thickness $2\epsilon^\kappa$)
tend to zero and positive constant $\beta-\alpha$, respectively.
\item[(ii)] We want to emphasize that
Theorem~\ref{NE-thm}(ii) improves the asymptotic behavior of $\phie(x)-\phie(\pm1)$ shown in Theorem 1.5 of our previous paper \cite{LHLL}.\\
\end{itemize}
\end{rem}

Following results play important roles throughout this paper.
\begin{itemize}
\item[(a)] Multiplying the equation (\ref{2-eqn1}) by $\phie'$ , (\ref{2-eqn1}) may be transformed into
\begin{eqnarray}
\frac{\epsilon^2}{2}\phie'^2(x)&=&\sumk\frac{\alpha_k}{\intt\eapy}\eapx
+\suml\frac{\beta_l}{\intt\ebpy}\ebpx\nonumber\\[-0.7em]
&&\label{2-eqn3}\\[-0.7em]
&&~~~+C_\epsilon,\nonumber
\end{eqnarray}
where $C_\epsilon$ is a constant depending on $\epsilon$.

\item[(b)] Differentiating (\ref{2-eqn1}) to $x$ and multiplying it by
$\phie'$,
\begin{eqnarray}
&&\epsilon^2\,\phie'''(x)\phie'(x)\nonumber\\[-0.7em]
&&\label{2-eqn4}\\[-0.7em]
&&=
\left(\sumk\frac{a_k^2\alpha_k}{\intt\eapy}\eapx
+\suml\frac{b_l^2\beta_l}{\intt\ebpy}\ebpx\right)\phie'^2(x).\nonumber
\end{eqnarray}
\end{itemize}

The rest of this paper is organized as follows: The proof of
Theorems \ref{2-thm5} and \ref{2-thm:v} are shown in
Section~\ref{2-mpnpsec:2}. In Section~\ref{2-mpnpsec:4}, we
compare the CCPB equation (\ref{2-eqn1}) and the PB equation
(\ref{2-eqnw}), and give the proof of Theorem~\ref{2-thm:w}.
In Section~\ref{NE-sec}, we consider the non-electroneutral case
and give the proof of Theorem~\ref{NE-thm}.
 In Section~\ref{2-mpnpsec:5}, several numerical experiments results
of the CCPB equation (\ref{2-eqn1}) and the PB equation
(\ref{2-eqnw}) are presented. The numerical computations are
basically preformed using finite element discretizations. In the
final section, we state the conclusion.

\section{Electroneutral cases: Proof of Theorems \ref{2-thm5} and \ref{2-thm:v} }\label{2-mpnpsec:2}
\ \ \ \ Let $\phie$ be the solution of the equation (\ref{2-eqn1})
with the boundary condition (\ref{2-eqn2}). A crucial property of $\phie$ is given as follows:\\
\begin{prop}\label{prop1.2}
Let $\phie\in
C^\infty((-1,1))\cap C^2([-1,1])$ be the solution of the equation
(\ref{2-eqn1}) with the boundary condition (\ref{2-eqn2}). Then the following properties hold.
\begin{itemize}
\item[(i)] Either $\phie'$ has at most one zero in $[-1,1]$, or $\phie'\equiv0$ on $[-1,1]$.
\item[(ii)] If $\phie$ is nontrivial (i.e., nonzero solution), then
\begin{equation}\label{id:2-005}
 \phie''(x_2) \phie'(x_2)>\phie''(x_1) \phie'(x_1) \quad \hbox{for}\quad -1 \leq x_1 < x_2 \leq 1.
\end{equation}
\end{itemize}
\end{prop}
\begin{proof}
We prove (i) by contradiction. Suppose there exist $y_1,\,y_2\in[-1,1]$ such that $y_1< y_2$ and $\phie'(y_1)=\phie'(y_2)=0$.
Then integrating (\ref{2-eqn4}) from $y_1$ to $y_2$ and using integration by parts, we may get
$$\int_{y_1}^{y_2}\,\epsilon^2\,(\phie'')^2\,+\left(\sumk\frac{a_k^2\alpha_k}{\intt\eapy}\eapx
+\suml\frac{b_l^2\beta_l}{\intt\ebpy}\ebpx\right)\phie'^2(x)\,dx=0\,,  $$ which implies $\phie'\equiv\phie''\equiv0$ on $[y_1,y_2]$. Here we have used the hypothesis $\phie'(y_1)=\phie'(y_2)=0$ and each $\alpha_k, \beta_l>0$.
On the other hand, the CCPB equation (\ref{2-eqn1}) has the following form $$\epsilon^2 \phie''(x) = \sumk A_{k,\epsilon} \,\eapx -\suml B_{l,\epsilon}\,\ebpx\quad\hbox{ for }\:x\in (-1,1)\,,$$ where $A_{k,\epsilon}$'s and $B_{l,\epsilon}$'s are constants, therefore $\phie$ satisfies the unique continuation property. Therefore, $\phie'$ has to be identically zero on $[-1,1]$. This completes the proof of Proposition~\ref{prop1.2}(i).

To prove (ii), we assume that $\phie$ is a nonzero solution of (\ref{2-eqn1}).
Thus, for any subinterval $(x_1,x_2)\subset(-1,1)$, Proposition~\ref{prop1.2}(i) immediately implies
\begin{equation}\label{2014-0808-1}
\int_{x_1}^{x_2}\left(\sumk\frac{a_k^2\alpha_k}{\intt\eapy}\eapx
+\suml\frac{b_l^2\beta_l}{\intt\ebpy}\ebpx\right)\phie'^2(x)dx>0.
\end{equation}
Integrating (\ref{2-eqn4}) over the interval $(x_1,x_2)$ and using (\ref{2014-0808-1}), we obtain (\ref{id:2-005})
and complete the proof of Proposition~\ref{prop1.2}.
\end{proof}
\medskip

The following interior estimate of $\phie$ is a key step for the proof of Theorem~\ref{2-thm5}.\\
\medskip
\noindent
\begin{theorem}\label{2-thm2}
Under the same hypotheses of Theorem~\ref{2-thm5}, we have
\begin{enumerate}
 \item[(i)]~~~ $-\phie(-1)=\phie(1)>0$ and $\phie'(1)=\phie'(-1)$.
 The solution $\phie$ is monotone increasing on $[-1,1]$, concave on $(-1,x_\epsilon^*)$ and convex on $(x_\epsilon^*,1)$,
 where $x_\epsilon^*\in(-1,1)$.
 Moreover, we have
 \begin{equation}
 \max_{x\in[-1,1]}|\phie(x)|\leq\phie(1)\leq\phi_0^+.\label{id:2010-1110}
 \end{equation}

 \item[(ii)]~~~ There are positive constants $C_1$ and $M_1$ independent of $\epsilon$ such that for
 $x \in [-1, 1]$ and $0<\epsilon\ll1$,
 \begin{equation}\label{id:2-006}
   0 \leq \phie'(x) \leq \frac{C_1}{\epsilon} \left( e^{-\frac{M_1(1+x)}{\epsilon}} +e^{-\frac{M_1(1-x)}{\epsilon}} \right).
  \end{equation}
\end{enumerate}
\end{theorem}

\begin{rem}\label{2-rk1.2}
\begin{itemize}
 \item[(i)] Replacing $\phie$ by $\phie+C$ for any constant
 $C$, the equation~(\ref{2-eqn1}) is invariant. Hence Theorem~\ref{2-thm2}~(i) implies that for any $\phi_0^+$ and
$\phi_0^-$, $\phie$ is monotonic on $[-1,1]$.
\item[(ii)] When
$N_1=N_2$, $\alpha_k=\beta_k$ and $a_k=b_k$ for $k=1,..., N_1$, as
for Theorem~1.2 in~\cite{LHLL}, the solution $\phie$ of
(\ref{2-eqn1}) and (\ref{2-eqn2}) is an odd function on $[-1,1]$,
and all denominator terms of (\ref{2-eqn1}) become equal. Then one
may follow the argument of~\cite{LHLL} to get the asymptotic
behavior of $\phie$'s. However, as $N_1\neq N_2$, $\alpha_k\neq
\beta_k$ or $a_k\neq b_k$ for some $k$, the solution $\phie$ may
not be odd on $[-1,1]$ so the argument of~\cite{LHLL} may fail for
this case and we have to develop a new argument to
prove~Theorem~\ref{2-thm2}.\\
\end{itemize}
\end{rem}

\subsection{Proof of Theorem \ref{2-thm2}}
\medskip
\noindent

Integrating (\ref{2-eqn1}) over $(-1,1)$ gives
$\intt\phie''(x)dx=\sumk\,a_k\alpha_k-\suml\,b_l\beta_l=0$. This
implies $\phie'(1)=\phie'(-1)$ and there exists
$x_\epsilon^*\in(-1,1)$ such that $\phie''(x_\epsilon^*)=0$. Along
with the boundary condition (\ref{2-eqn2}), we find
$-\phie(1)=\phie(-1)$. Setting $x_1=x_\epsilon^*$ and
$x_2=x_\epsilon^*$ in (\ref{id:2-005}), respectively, we get
\begin{equation}\label{1.7.1}
\phie''(x)\phie'(x)>0\mbox{ for }
x\in(x_\epsilon^*,1],\,\mbox{ and }
\phie''(x)\phie'(x)<0 \mbox{ for }
x\in[-1,x_\epsilon^*),
\end{equation}
which implies: \textbf{(a)} Both $\phie'$ and $\phie''$ never change sign and
 share the same sign on $(x_\epsilon^*,1]$;
\textbf{(b)}~$\phie'$ and $\phie''$ never change sign and have opposite signs on
$[-1,x_\epsilon^*)$. Consequently,
$\phie'$ never changes sign on $[-1,1]$ due to (a), (b) and $\phie'(1)=\phie'(-1)$.
Now we claim $\phie'\geq0$ on $[-1,1]$. We state the proof using contradiction.
Suppose $\phie'\leq0$ on $[-1,1]$, then the boundary condition~(\ref{2-eqn2}) implies
$\phi_0^-=\phie(-1)-\eta_\epsilon\phie'(-1)\geq\phie(1)+\eta_\epsilon\phie'(1)=\phi_0^+$,
which gives a contradiction.

Therefore, we get $\phie'\geq0$ on $[-1,1]$.
Along with the boundary condition~(\ref{2-eqn2}), we prove (\ref{id:2010-1110}).


Furthermore, by
(\ref{1.7.1}) and $\phie'(x)\geq 0$, we have
 $\phie''(x)\leq0$ for $x\in(-1,x_\epsilon^*)$ and
$\phie''(x)\geq0$ for $x\in(x_\epsilon^*,1)$.
Hence we complete the proof of Theorem \ref{2-thm2} (i).

By (\ref{id:2010-1110}) and (\ref{2-eqn4}), we obtain
\begin{equation}\label{1.8.1}
\epsilon^2(\phie'^2(x))''
\geq2\epsilon^2\phie'''(x)\phie'(x)\geq\,4M_1^2\phie'^2(x)
\end{equation}
for $x\in(-1,1)$ and $\epsilon>0$, where
$M_1\!=\!\frac{1}{2}\left(\sumk\,a_k^2\alpha_ke^{2a_k\phi_0^-}
\!+\!\suml\,b_l^2\beta_le^{2b_l\phi_0^-}\right)^{1/2}$.\\
Note that $\phie'(1)=\phie'(-1)>0$. By~(\ref{1.8.1}) and the
standard comparison theorem, we get
 \begin{equation}\label{id:2-2010-1110-1}
 0\leq\phie'(x) \leq \phie'(1)\left( e^{-\frac{M_1(1+x)}{\epsilon}} +e^{-\frac{M_1(1-x)}{\epsilon}} \right).
  \end{equation}
It remains to deal with $\phie'(1)$. By (\ref{id:2010-1110}),
there exists $\overline{x}_\epsilon\in(-1,1)$ such that
$0\leq\phie'(\overline{x}_\epsilon)=\frac{\phie(1)-\phie(-1)}{2}\leq\phi_0^+$.
Subtracting (\ref{2-eqn3}) at $x=\overline{x}_\epsilon$ from that
at $x=1$ and using (\ref{id:2010-1110}), it is easy to get
$\phie'(1)\leq\frac{C_1}{\epsilon}$ as $0<\epsilon\ll1$, where
$C_1$ is a positive constant independent of $\epsilon$. Along with
(\ref{id:2-2010-1110-1}), we get (\ref{id:2-006}) and prove
Theorem \ref{2-thm2}(ii).

Therefore, we complete the proof of Theorem \ref{2-thm2}.\\

Note that (\ref{2-eqn3}) plays a crucial role on the asymptotic
behavior of $\phie$ as $\epsilon\downarrow0$. The estimate of the
constant $C_\epsilon$ in (\ref{2-eqn3}) is given as follows:\\
\begin{lem}\label{2-lem3}
Under the same hypotheses of Theorem~\ref{2-thm2}, we have
\begin{enumerate}
\item[(i)]~~~For any $x,y\in(-1,1)$, $\phie(x)-\phie(y)$ converges
exponentially to zero as $\epsilon$ goes to zero.
\item[(ii)]~~$\ds\lim_{\epsilon\downarrow0}\, C_\epsilon=
-\frac{1}{2}\left(\sumk\alpha_k +\suml\beta_l\right)$, where
$C_\epsilon$ is the constant defined in~(\ref{2-eqn3}).
\end{enumerate}
\end{lem}
\begin{proof}
(\ref{id:2-006}) implies
that for any $x,\,y\in(-1,1)$,
\begin{eqnarray}
&&\lim_{\epsilon\downarrow0}\phie'(x)= 0\,\mbox{ and }\,\nonumber\\[-0.7em]
&&\label{id:2-007}\\[-0.7em]
&&|\phie(x)-\phie(y)|\leq\frac{C_1}{M_1}\!
\left(\left|e^{-\!\frac{M_1(1+x)}{\epsilon}}\!+\!e^{-\!\frac{M_1(1-x)}{\epsilon}}\right|
\!+\!\left|e^{-\!\frac{M_1(1+y)}{\epsilon}}\!+\!e^{-\!\frac{M_1(1-y)}{\epsilon}}\right|\right).\nonumber
\end{eqnarray}
This may complete the proof of Lemma~\ref{2-lem3}(i). Note that
(\ref{id:2-007}) gives
$\sup_{x,\,y\in(-1,1)}|\phie(x)-\phie(y)|\leq4C_1/M_1$ and
$\lim_{\epsilon\downarrow0}|\phie(x)-\phie(y)|=0$ for
$x,\,y\in(-1,1)$. Applying Lebesgue's dominated convergence
theorem, we obtain
\begin{equation}\label{id:2-009}
\lim_{\epsilon\downarrow0}\,\frac{\alpha_k\eapx}{\intt\eapy}
=\frac{\alpha_k}{2}\, \quad\mathrm{and}\quad
\lim_{\epsilon\downarrow0}\,\frac{\beta_l\ebpx}{\intt\ebpy}
=\frac{\beta_l}{2}\,,
\end{equation}
for $k=1,\cdots,N_1$ and $l=1,\cdots,N_2$.

Therefore, by~(\ref{2-eqn3}), (\ref{id:2-007}) and (\ref{id:2-009}), we prove
Lemma \ref{2-lem3}(ii) and complete the proof of~Lemma
\ref{2-lem3}.
\end{proof}

\subsection{Proof of Theorem \ref{2-thm5}}
To prove Theorem~\ref{2-thm5}, we need the following lemma:\\
\begin{lem}\label{2-lem4}
(i) Under the same hypotheses of Theorem~\ref{2-thm2}, we have

\begin{equation}\label{id:2-011}
\sumk\frac{\alpha_k\left(e^{a_k\phie(1)}-e^{-a_k\phie(1)}\right)}{\intt\eapy}
=\suml\frac{\beta_l\left(e^{b_l\phie(1)}-e^{-b_l\phie(1)}\right)}{\intt\ebpy}.
\end{equation}
$(ii)$ If $\eta_\epsilon\neq 0$, then
\begin{equation}\label{id:2-012}
\frac{\epsilon^2}{2\eta_\epsilon^2}(\phi_0^+-\phie(1))^2
=\sumk\frac{\alpha_ke^{a_k\phie(1)}}{\intt\eapy}+\suml\frac{\beta_le^{-b_l\phie(1)}}{\intt\ebpy}
+C_\epsilon.
\end{equation}
\end{lem}

\begin{proof}
To get (\ref{id:2-011}), we subtract the equation
(\ref{2-eqn3}) at $x=-1$ from that at $x=1$. Here we have used the
facts that $\phie'(1)=\phie'(-1)$ and $\phie(-1)=-\phie(1)$ which
come from Theorem~\ref{2-thm2}~(i). Setting $x=1$ in
(\ref{2-eqn3}), we use (\ref{2-eqn2}) to get (\ref{id:2-012}),
and complete the proof of Lemma~\ref{2-lem4}.
\end{proof}

To uniquely determine the values $c$ and $t$, we need the
following lemma.\\
\begin{lem}\label{2-lem7}
Assume $\sumk\,a_k\alpha_k=\suml\,b_l\beta_l$. Then
\begin{enumerate}
\item[(i)]~~~$f$ is strictly increasing on $(0,\infty)$ and
strictly decreasing on $(-\infty,0)$.
\item[(ii)]~There exists a
unique solution $(t,c)$ of the equations
(\ref{1.19.1})-(\ref{1.19.3}).
\end{enumerate}
\end{lem}
\begin{proof}
By (\ref{id:2-015}) and $\sumk\,a_k\alpha_k=\suml\,b_l\beta_l$, it
is easy to check that $f'(s)>0$ if $s>0$ and $f'(s)<0$ if $s<0$.
This shows~(i). To prove (ii), we need

\textbf{Claim 1.} There exists
$0<s<2\phi_0^+-2\gamma(f(s)-f(0))^{1/2}$ such that
$$f(s)=f\left(s-2\phi_0^++2\gamma(f(s)-f(0))^{1/2}\right)\,.$$
\begin{proof}[Proof of Claim 1]
Let $k(s)=s-2\phi_0^++2\gamma(f(s)-f(0))^{1/2}$ for
$s\in\mathbb{R}$. Then $k(0)=-2\phi_0^+<0$, $k(2\phi_0^+)>0$ and
$k'(s)=1+\gamma(f(s)-f(0))^{-1/2}\,f'(s)>0$ for $s\in (0,\infty)$.
Hence, there exists $s_1\in(0,2\phi_0^+)$ such that $k(s_1)=0$
and
\begin{equation}\label{1.19.7}
k(s)<0\quad\hbox{ for }\quad
s\in(0,s_1)\,.
\end{equation}
Let $h(s)=f(s)-f(k(s))$ for $s\in\mathbb{R}$. Then
$h(0)=f(0)-f(-2\phi_0^+)<0$ and $h(s_1)=f(s_1)-f(0)>0$. Hence
there exists $s_2\in(0,s_1)$ such that $h(s_2)=0$. On the other
hand, (\ref{1.19.7}) implies $k(s_2)<0$
i.e.~$s_2<2\phi_0^+-2\gamma(f(s_2)-f(0))^{1/2}$.

Therefore, we complete the proof of Claim~1.
\end{proof}

\medskip
Now we want to prove Lemma~\ref{2-lem7}(ii). By Claim~1,
$$(t,c)=\left(\phi_0^+-\gamma(f(s)-f(0))^{1/2},\phi_0^+-s-\gamma(f(s)-f(0))^{1/2}\right)$$
is a solution of (\ref{1.19.1}) and (\ref{1.19.2}). Moreover,
$0<s<2\phi_0^+-2\gamma(f(s)-f(0))^{1/2}$ gives $|c|<t\leq\phi_0^+$.
Hence (\ref{1.19.1})-(\ref{1.19.3}) have a solution. The
uniqueness of~(ii) can be proved by contradiction. Suppose
$(t_1,c_1)$ and $(t_2,c_2)$ solve (\ref{1.19.1})-(\ref{1.19.3})
and $t_1>t_2$. If $f(t_1-c_1)\geq\,f(t_2-c_2)$, then
$\phi_0^+-t_2>\phi_0^+-t_1=\gamma(f(t_1-c_1)-f(0))^{1/2}\geq\gamma(f(t_2-c_2)-f(0))^{1/2}$
i.e.~$\phi_0^+-t_2>\gamma(f(t_2-c_2)-f(0))^{1/2}$ contradicts
$(t_2,c_2)$ a solution of (\ref{1.19.1})-(\ref{1.19.3}). Thus
$f(t_1-c_1)<f(t_2-c_2)$ and then (\ref{1.19.2}) gives
$f(-t_1-c_1)<f(-t_2-c_2)$. Furthermore, by Lemma~\ref{2-lem7}(i),
we obtain $t_1-c_1<t_2-c_2$ and $-t_1-c_1>-t_2-c_2$ which implies
$t_1<t_2$ a contradiction to the hypothesis $t_1>t_2$. Hence, $t_1=t_2:=t^*$. Here we
have used the facts that $t_1-c_1,t_2-c_2>0$ and
$-t_1-c_1,-t_2-c_2<0$.

To prove $c_1=c_2$, we set $g(s):=f(t^*-s)-f(-t^*-s)$. Note that
Lemma~\ref{2-lem7}~(i) implies $g'(s)=-f'(t^*-s)+f'(-t^*-s)<0$ for
$|s|<t^*$, i.e., $g(s)$ is strictly decreasing on $(-t^*,t^*)$.
Therefore, we have $c_1=c_2$ and complete the proof of
Lemma~\ref{2-lem7}~(ii).
\end{proof}
\medskip

Now we shall give the proof of Theorem~\ref{2-thm5}.\\

\textbf{Proof of Theorem~\ref{2-thm5}.}
 By Lemma~\ref{2-lem3}(i), it suffices to prove
$\ds\lim_{\epsilon\downarrow0}\phie(0)=c$. By
Theorem~\ref{2-thm2}, $\left\{|\phie(0)|\right\}_{\epsilon>0}$ has
an upper bound. Then we set
$\displaystyle\limsup_{\epsilon\downarrow0}\phie(0)=c_s$ and
$\displaystyle\liminf_{\epsilon\downarrow0}\phie(0)=c_i$. Hence
there exist sequences $\{\epsilon_j\}_{j\in \mathbb{N}}$ and
$\{\widetilde{\epsilon}_{j}\}_{j\in \mathbb{N}}$ tending to zero
such that
$\displaystyle\lim_{j\rightarrow\infty}\phi_{\epsilon_j}(0)=c_s$
and
$\displaystyle\lim_{j\rightarrow\infty}\phi_{\widetilde{\epsilon}_j}(0)=c_i$.
We may rewrite (\ref{id:2-011}) and (\ref{id:2-012}) as follows:
\begin{eqnarray}
&&\sumk\frac{\alpha_k\left(e^{a_k(\phie(1)-\phie(0))}-e^{-a_k(\phie(1)+\phie(0))}\right)}{\intt\,e^{a_k(\phie(y)-\phie(0))}dy}\nonumber\\[-0.7em]
&&\label{id:2-013}\\[-0.7em]
&&~~~~~~=\suml\frac{\beta_l\left(e^{b_l(\phie(1)+\phie(0))}-e^{-b_l(\phie(1)-\phie(0))}\right)}{\intt\,e^{-b_l(\phie(y)-\phie(0))}dy}\nonumber
\end{eqnarray}
and
\begin{eqnarray}
\frac{\epsilon^2}{2\eta_\epsilon^2}(\phi_0^+-\phie(1))^2
&=&\sumk\frac{\alpha_ke^{a_k(\phie(1)-\phie(0))}}{\intt\,e^{a_k(\phie(y)-\phie(0))}dy}\nonumber\\[-0.7em]
&&\label{id:2-014}\\[-0.7em]
&&~~~+\suml\frac{\beta_le^{-b_l(\phie(1)-\phie(0))}}{\intt\,e^{-b_l(\phie(y)-\phie(0))}dy}\nonumber
+C_\epsilon.
\end{eqnarray}
We divide the proof into two cases.\\

\textbf{Case 1}. $\lim_{\epsilon\downarrow0}\frac{\eta_\epsilon}{\epsilon}=\infty$.\\
Note that $|\phi_0^+-\phie(1)|\leq2\phi_0^+$. By
(\ref{id:2-009}), (\ref{id:2-015}), (\ref{id:2-014}) and
Lemma~\ref{2-lem3}~(ii), we have
$f(\limsup_{\epsilon\downarrow0}(\phie(1)-\phie(0)))=f(\liminf_{\epsilon\downarrow0}(\phie(1)-\phie(0)))=f(0)$.
Then by Lemma~\ref{2-lem7}~(i),
$\lim_{\epsilon\downarrow0}(\phie(1)-\phie(0))=0$. Along with
(\ref{id:2-013}), we find
$f(-\lim_{\epsilon\downarrow0}(\phie(1)+\phie(0)))=f(0)$, this
gives $\lim_{\epsilon\downarrow0}(\phie(1)+\phie(0))=0$.
Consequently, we have
$\lim_{\epsilon\downarrow0}\phie(1)=\lim_{\epsilon\downarrow0}\phie(0)=0$.
Hence, we obtain $c=t=0$ and complete the proof of
Theorem~\ref{2-thm5}~(i).
\\

\textbf{Case 2}.
$\lim_{\epsilon\downarrow0}\frac{\eta_\epsilon}{\epsilon}=\gamma<\infty$.\\
By Theorem~\ref{2-thm2},
$\left\{|\phi_{\epsilon_j}(1)|\right\}_{j\in \mathbb{N}}$ has an
upper bound. Then there is a constant $t_s$ and a subsequence of
$\{\epsilon_j\}$ (for notation convenience, we still denote it by
$\{\epsilon_j\}$) such that
$\lim_{j\rightarrow\infty}\phi_{\epsilon_j}(1)=t_s$. Putting
$\epsilon=\epsilon_j$ in (\ref{id:2-013}) and (\ref{id:2-014}) and using Lemma~\ref{2-lem3}(ii), one may check that $(t_s,c_s)$ satisfies
\begin{equation}\label{id:2-019}
(\phi_0^+-t_s)^2=\gamma^2(f(t_s-c_s)-f(0))\quad\mathrm{and}\quad
f(t_s-c_s)=f(-t_s-c_s).
\end{equation}

Now we claim $|c_s|<t_s\leq\phi_0^+$.
Since
$|\phie(0)|\leq\phi(1)\leq\phi_0^+\neq0$, then we have
$|c_s|\leq t_s\leq\phi_0^+\neq0$. If $|c_s|=t_s$, then by
the second equation of (\ref{id:2-019}) and Lemma~\ref{2-lem7}(i), we have
$t_s-c_s=-t_s-c_s=0$, i.e. $t_s=c_s=0$. Along with the first equation of (\ref{id:2-019})
we find $\phi_0^+=0$, which contradicts to
$\phi_0^+\neq0$. Hence $|c_s|<t_s\leq\phi_0^+$.
Along
with (\ref{id:2-019}), $(t_s,c_s)$ satisfies
(\ref{1.19.1})-(\ref{1.19.3}).

Similarly, there is a positive
constant $t_i$ such that $(t_i,c_i)$ satisfies
(\ref{1.19.1})-(\ref{1.19.3}). By Lemma~\ref{2-lem7}~(ii), we get
$c_s=c_i=c$ and $t_s=t_i=t$, where
$\lim_{\epsilon\downarrow0}\phie(0)=c$,
$\lim_{\epsilon\downarrow0}\phie(1)=t$ and $(t,c)$ satisfies
(\ref{1.19.1})-(\ref{1.19.3}). Therefore, we may complete the
proof of Theorem \ref{2-thm5}(ii).

By (\ref{id:2-015}) and $|t-c|\leq2\phi_0^+$, $f(t-c)-f(0)$ is uniformly bounded for all $\gamma>0$.
Consequently, by (\ref{1.19.1})-(\ref{1.19.3}), we have $\lim_{\gamma\rightarrow0}t(\gamma)=\phi_0^+$ and $\lim_{\gamma\rightarrow0}c(\gamma)=c^*$,
where $|c^*|<\phi_0^+$ is uniquely determined by $f(\phi_0^+-c^*)= f(-\phi_0^+-c^*)$.
By (\ref{1.19.1}) we have $f(t-c)-f(0)=\left(\frac{\phi_0^+-t}{\gamma}\right)^2$, which and (\ref{1.19.2}) give $\lim_{\gamma\rightarrow\infty}f(t-c)=\lim_{\gamma\rightarrow\infty}f(-t-c)=f(0)$.
By Lemma \ref{2-lem7}(i) and the continuity of $f$, we find $\lim_{\gamma\rightarrow\infty}(t-c)=\lim_{\gamma\rightarrow\infty}(-t-c)=0$.
Hence, $\lim_{\gamma\rightarrow\infty}t=\lim_{\gamma\rightarrow\infty}c=0$
and complete the proof of Theorem~\ref{2-thm5}(ii-A).

It remains to prove Theorem \ref{2-thm5}(ii-B). By
(\ref{1.19.1})-(\ref{1.19.3}), $t$ and $c$ are uniquely determined
by $\gamma$. Hence we can consider $t=t(\gamma)$ and $c=c(\gamma)$
as functions of $\gamma$. Due to
$f(t(\gamma)-c(\gamma))-f(0)\neq0$ on $(0,\infty)$ (by
(\ref{1.19.3}) and Lemma \ref{2-lem7}(i)),
$t(\gamma),\,c(\gamma):\,
(0,\infty)\rightarrow(-\phi_0^+,\phi_0^+)$ are continuously
differentiable. Differentiating (\ref{1.19.1}) and (\ref{1.19.2})
to $\gamma$, one may check that
\begin{equation}\label{check1}
\left[1-\gamma(f(t-c)-f(0))^{-1/2}\frac{f'^2(t-c)}{f'(t-c)-f'(-t-c)}\right]\frac{dt}{d\gamma}=-(f(t-c)-f(0))^{1/2}
\end{equation}
and
\begin{equation}\label{check2}
\frac{d}{d\gamma}(t-c)= -\frac{2f'(-t-c)}{f'(t-c)-f'(-t-c)}\frac{dt}{d\gamma}.
\end{equation}
If $\frac{dt}{d\gamma}$ changes the sign on $(0,\infty)$, then there is
a $\gamma^*\in(0,\infty)$ such that
$\frac{dt}{d\gamma}(\gamma^*)=0$. By (\ref{check1}) and Lemma
\ref{2-lem7}(i), we have $t(\gamma^*)=c(\gamma^*)$ contradicting
to (\ref{1.19.3}). Hence $\frac{dt}{d\gamma}$ keeps the same sign
on $(0,\infty)$, and then Theorem ~\ref{2-thm5}(iii) gives
$\frac{dt}{d\gamma}<0$ on $(0,\infty)$. On the other hand,
Lemma~\ref{2-lem7}(i) and~(\ref{1.19.3}) imply that both
$f'(t-c)$ and $-f'(-t-c)$ are positive. Consequently, by
(\ref{check2}), we obtain that $\frac{d}{d\gamma}(t-c)$ and
$\frac{dt}{d\gamma}$ share the same sign.
Therefore, we prove Theorem~\ref{2-thm5}~(ii-B) and complete the proof of Theorem~\ref{2-thm5}.\\
\begin{rem}\label{2-rk4}
Suppose
$\lim_{\epsilon\downarrow0}\frac{\epsilon}{\eta_\epsilon}=0$. Then
Theorem \ref{2-thm2}~(i) and Theorem~\ref{2-thm5}~(i) give
$\phie\rightarrow 0$ uniformly in $[-1,1]$ as
$\epsilon\downarrow0$.
\end{rem}

\subsection{Proof of Theorem~\ref{2-thm:v}}
\medskip
\noindent

For convenience, setting $v(x)=\phie(x)-c$,  By (\ref{2-eqn3}), we
find
\begin{equation}\label{2-v1}
\frac{\epsilon^2}{2}v'^2(x)=\sumk\frac{\alpha_k}{\intt\,e^{a_kv(y)}dy}e^{a_kv(x)}
+\suml\frac{\beta_l}{\intt\,e^{-b_lv(y)}dy}e^{-b_lv(x)}
+C_\epsilon.
\end{equation}
Note that by Theorem \ref{2-thm2}(i) and Theorem \ref{2-thm5}, we
have $|v(x)|\leq\phi_0^++|c|$ and
$\ds\lim_{\epsilon\downarrow0}\intt\,e^{a_kv(y)}dy=
\lim_{\epsilon\downarrow0}\intt\,e^{-b_lv(y)}dy=2$. Hence by
Lemma~\ref{2-lem3}(ii), it is easy to check that
\begin{equation}\label{2-v4}
\left|\epsilon^2v'^2(x)-\left[f(v(x))-f(0)\right]\right|\leq\delta(\epsilon),
\end{equation}
for all $x\in[-1,1]$, where $\delta(\epsilon)$ is a positive
quantity tending to zero as $\epsilon$ goes to zero. By
Theorem~\ref{2-thm2}(i) and Theorem~\ref{2-thm5}(ii),
we have $v'=\phie'\geq0$ on $(-1,1)$ and
\begin{equation}\label{2-vtc}
\lim_{\epsilon\downarrow0}v(1)=t-c>0>-t-c=\lim_{\epsilon\downarrow0}v(-1).
\end{equation}
Thus there exist $\epsilon^*>0$ and $y_\epsilon^+\in(-1,1)$  such
that
$v(y_\epsilon^+)=-\log\left\{1-\left[\delta(\epsilon)/(\alpha_1+\beta_2)^{2}\right]^{1/4}\right\}>0$
for $0<\epsilon<\epsilon^*$. Note that
$\lim_{\epsilon\downarrow0}v(y_\epsilon^+)=0$ and
\begin{equation}\label{2-v5+}
v(x)\geq
v(y_\epsilon^+)=-\log\left\{1-\left[\delta(\epsilon)/(\alpha_1+\beta_2)^{2}\right]^{1/4}\right\}>0,\quad\forall\,
x\in(y_\epsilon^+,1).
\end{equation}

Now we begin to deal with (\ref{2-v4}) when $N_1=1$, $N_2=2$,
$a_1=b_1=1$ and $b_2=2$. Here we have $\alpha_1=\beta_1+2\beta_2$
and $f(v(x))-f(0)=(1-e^{-v(x)})^2(\alpha_1e^{v(x)}+\beta_2)$. Note
that such a formula is valid only when $N_1=1$, $N_2=2$,
$a_1=b_1=1$ and $b_2=2$. Along with (\ref{2-v5+}), we have
\begin{equation}\label{2-v_20101105}
K(\epsilon)(1\!-\!e^{-v(x)})^2(\alpha_1e^{v(x)}\!+\!\beta_2)\!+\!\delta(\epsilon)\leq\!f(v(x))\!-\!f(0)
\leq\!(t_\epsilon\!-\!e^{-v(x)})^2(\alpha_1e^{v(x)}\!+\!\beta_2)\!-\!\delta(\epsilon),
\end{equation}
for $x\in(y_\epsilon^+,1)$ and $0<\epsilon<\epsilon^*$, where
$K(\epsilon)=1-\sqrt{\delta(\epsilon)}$ and
$t_\epsilon=1+\sqrt{\delta(\epsilon)/(\alpha_1+\beta_2)}$.
Consequently, (\ref{2-v4})--(\ref{2-v_20101105}) give
\begin{equation}\label{2-v9}
\frac{v'(x)}{(t_\epsilon-e^{-v(x)})\sqrt{\alpha_1e^{v(x)}+\beta_2}}\leq\frac{1}{\epsilon}
\leq\frac{1}{\sqrt{K(\epsilon)}}\cdot\frac{v'(x)}{(1-e^{-v(x)})\sqrt{\alpha_1e^{v(x)}+\beta_2}},
\end{equation}
for $x\in(y_\epsilon^+,1)$.

Integrate (\ref{2-v9}) over $(y,1)$ for
$y\in(y_\epsilon^+,1)$, we obtain
\begin{eqnarray}
&&\int_y^1\frac{v'(x)}{(a-e^{-v(x)})\sqrt{\alpha_1e^{v(x)}+\beta_2}}dx\nonumber\\[-0.7em]
&&\label{20110518}\\[-0.7em]
&=&
\frac{1}{\sqrt{\alpha_1a+\beta_2a^2}}
\log\frac
{(\sqrt{(\alpha_1e^{v(1)}+\beta_2)a}-\sqrt{\alpha_1+\beta_2a})(\sqrt{(\alpha_1e^{v(y)}+\beta_2)a}+\sqrt{\alpha_1+\beta_2a})}
{(\sqrt{(\alpha_1e^{v(1)}+\beta_2)a}+\sqrt{\alpha_1+\beta_2a})(\sqrt{(\alpha_1e^{v(y)}+\beta_2)a}-\sqrt{\alpha_1+\beta_2a})},\nonumber
\end{eqnarray}
for $a>0$. Hence, (\ref{2-v9}) and (\ref{20110518}) imply
\begin{eqnarray}
&&t_\epsilon\!+\!\!\left(t_\epsilon\!+\!\frac{\beta_2}{\alpha_1}\right)\!
\mathrm{csch}^{2}\!\!\left[\!\frac{C_{1,\epsilon}^+}{\epsilon}(1\!-\!x)\!+\!\log\,D_{1,\epsilon}^+\!\right]\!
\leq\!e^{v(x)}\nonumber\\[-0.7em]
&&\label{2-v11}\\[-0.7em]
&&~~~~~~~~~~~~~~~~~~~~\leq\!1\!+\!\left(1\!+\!\frac{\beta_2}{\alpha_1}\right)\!
\mathrm{csch}^{2}\!\!\left[\!\frac{C_{2,\epsilon}^+}{\epsilon}(1\!-\!x)\!+\!\log\,D_{2,\epsilon}^+\!\right],\nonumber
\end{eqnarray}
for $x\in(y_\epsilon^+,1)$, where
$C_{1,\epsilon}^+=\sqrt{\alpha_1t_\epsilon+\beta_2}$,
$C_{2,\epsilon}^+=K(\epsilon)\sqrt{\alpha_1+\beta_2}$,
$D_{1,\epsilon}^+=\frac{\sqrt{\alpha_1e^{v(1)}+\beta_2}+\sqrt{\alpha_1t_\epsilon+\beta_2}}
{\sqrt{\alpha_1e^{v(1)}+\beta_2}-\sqrt{\alpha_1t_\epsilon+\beta_2}}$
and
$D_{2,\epsilon}^+=\frac{\sqrt{\alpha_1e^{v(1)}+\beta_2}+\sqrt{\alpha_1+\beta_2}}
{\sqrt{\alpha_1e^{v(1)}+\beta_2}-\sqrt{\alpha_1+\beta_2}}$. By
(\ref{2-vtc}), (\ref{2-v11}) and
$\lim_{\epsilon\downarrow0}t_\epsilon=\lim_{\epsilon\downarrow0}K(\epsilon)=1$,
we get (\ref{2-id:thmw1}).

Similarly, we also have (\ref{2-id:thmv}).

Therefore, we complete the proof of
Theorem~\ref{2-thm:v}.

When $N_1=N_2=2$, $a_i=b_i=i$, $i=1,2$, and
$\alpha_1+2\alpha_2=\beta_1+2\beta_2$, we may follow the similar
proof of Theorem~\ref{2-thm:v} and obtain the following result.\\

\begin{cor}\label{2-cor2}
Under the same hypotheses of Theorem~\ref{2-thm2}, suppose
$N_1=N_2=2$, $a_i=b_i=i$, $i=1,2$, and
$\alpha_1+2\alpha_2=\beta_1+2\beta_2$. Then
\begin{eqnarray}
\phi_{1,\epsilon}^{+}(x)&\leq\,\phie(x)\leq\phi_{2,\epsilon}^{+}(x),\quad \forall\, x\in(\overline{x}_\epsilon,1),\label{2-id:ww8}\\
\phi_{1,\epsilon}^{-}(x)&\leq\,\phie(x)\leq\phi_{2,\epsilon}^{-}(x),\quad
\forall\, x\in(-1,\overline{x}_\epsilon)\,,\label{2-id:ww9}
\end{eqnarray}
where
$$
\phi_{i,\epsilon}^{\pm}(x)=c+\log\frac{\cosh
h_{i,\epsilon}^{\pm}(x)\pm\frac{A^2-B^2+A}{B}} {\cosh
h_{i,\epsilon}^{\pm}(x)\mp\frac{A+1}{B}},\quad
h_{i,\epsilon}^{\pm}(x)={\frac{\widetilde{C}_{i,\epsilon}^{\pm}}{\epsilon}(1\mp\,x)}+\log
H_{i,\epsilon}^{\pm},\, i=1,2\,.
$$
Here $A=1+\frac{\alpha_1}{2\alpha_2}$,
$B=\sqrt{\left(1+\frac{\alpha_1}{2\alpha_2}\right)^2-\frac{\beta_2}{\alpha_2}}$,
and $\widetilde{C}_{i,\epsilon}^\pm$'s, $H_{i,\epsilon}^{\pm}$'s,
$i=1,2$, are positive constants depending on $\epsilon$ such that
\begin{eqnarray}
\lim_{\epsilon\downarrow0}\widetilde{C}_{i,\epsilon}^\pm&=&\sqrt{\alpha_2[(A+1)^2-B^2]}\,,\nonumber\\
\lim_{\epsilon\downarrow0}H_{i,\epsilon}^{\pm}&=&
\left(\sqrt{\frac{A-B+e^{\pm\,t-c}}{A+B+e^{\pm\,t-c}}}+\sqrt{\frac{A-B+1}{A+B+1}}\right)\left(\pm\sqrt{\frac{A-B+e^{\pm\,t-c}}{A+B+e^{\pm\,t-c}}}
\mp\sqrt{\frac{A-B+1}{A+B+1}}\right)^{-1}.\nonumber
\end{eqnarray}
\end{cor}

\section{Proof of Theorem \ref{2-thm:w}}\label{2-mpnpsec:4}
\medskip
\ \ \ \ In this section, we study the asymptotic behavior of solution $\phie$ of the PB equation (\ref{2-eqnw}) with the boundary condition~(\ref{2-eqn2})
 and give the proof of Theorem \ref{2-thm:w}.
Surely, the PB equation~(\ref{2-eqnw}) can be transformed into
\begin{equation}\label{2-id:ww1}
\epsilon^2\phie''(x)=\frac{1}{2}f'(\phie(x))\,,
\end{equation}
where
$f(s)=\sum_{k=1}^{N_1}\alpha_ke^{a_ks}+\sum_{l=1}^{N_2}\beta_le^{-b_ls}$
is defined by (\ref{id:2-015}).
It is well-known that the equation~(\ref{2-eqnw})
has the unique solution $\phie\in\,C^\infty((-1,1))\cap\,C^2([-1,1])$.
As for (\ref{2-eqn3}), we use~(\ref{2-id:ww1}) to derive the
following identity
\begin{equation}\label{2-id:w2}
\frac{\epsilon^2}{2}\phie'^2(x)
=\frac{1}{2}f(\phie(x))+C'_\epsilon.
\end{equation}
Moreover, we use the similar argument of
(\ref{2-eqn4})-(\ref{id:2-005}) to get
\begin{equation}\label{id:2-314}
\phie''(x_2)\phie'(x_2)>\phie''(x_1)\phie'(x_1)\quad\mathrm{for}-1<x_1<x_2<1.
\end{equation}
Using standard maximum principle to
(\ref{2-eqnw}) and (\ref{2-eqn2}), we obtain
\begin{equation}\label{id:2-302}
\min\{0,\phi_0^+,\phi_0^-\}\leq\phie(x)\leq\max\{0,\phi_0^+,\phi_0^-\},
\end{equation}
for $x\in[-1,1]$.

Now we state the proof of Theorem~\ref{2-thm:w}.\\

\textbf{Proof of Theorem~\ref{2-thm:w}.}
 Multiplying the
equation (\ref{2-eqnw}) by $\phie$, we obtain
\begin{equation}\label{id:2-303}
\epsilon^2(\phie^2(x))''\geq 2\epsilon^2\phie''(x)\phie(x)
=\left(\sum_{k=1}^{N_1}a_k\alpha_ke^{a_k\phie(x)}-\sum_{l=1}^{N_2}b_l\beta_le^{-b_l\phie(x)}\right)\phie(x)
\geq C_5\,\phie^2(x)\,
\end{equation}
where $C_5=\displaystyle\inf_{s\in\mathbb{R}, s\neq
0}\,s^{-1}\left(\sum_{k=1}^{N_1}a_k\alpha_ke^{a_k\,s}-\sum_{l=1}^{N_2}b_l\beta_le^{-b_l\,s}\right)>0$.
Here we have used the hypothesis
$\sum_{k=1}^{N_1}a_k\alpha_k=\sum_{l=1}^{N_2}b_l\beta_l$ to assure
$C_5$ as a positive constant. Thus by (\ref{id:2-302}),
(\ref{id:2-303}) and the standard comparison theorem, we get
\begin{equation}\label{id:2-w}
|\phie(x)|\leq\max\{|\phi_0^+|,|\phi_0^-|\}\left(e^{-\frac{\sqrt{C_5}}{2\epsilon}(1+x)}+e^{-\frac{\sqrt{C_5}}{2\epsilon}(1-x)}\right)\,,
\quad\forall\:x\in (-1,1)\,,
\end{equation}
and complete the proof of Theorem~\ref{2-thm:w}(i).

Suppose $\min\{\phi_0^+,\phi_0^-\}>0$. Then (\ref{id:2-302})
gives $0\leq\phie(x)\leq\max\{\phi_0^+,\phi_0^-\}$, together
with (\ref{2-id:ww1}) and Lemma~\ref{2-lem7} (i), we may find
$\phie''\geq0$ on $[-1,1]$. Here we have used the hypothesis that
$\sum_{k=1}^{N_1}a_k\alpha_k=\sum_{l=1}^{N_2}b_l\beta_l$.
To complete the proof of Theorem~\ref{2-thm:w} (i), we need to claim:\\

\textbf{Claim 2.} There exist $\epsilon^*>0$ and
$x_\epsilon^*\in(-1,1)$ such that
$\ds\phie(x_\epsilon^*)=\min_{x\in[-1,1]}\phie(x)$ for
$0<\epsilon<\epsilon^*$.
\begin{proof}
We state the proof of Claim 2 by contradiction.
Suppose $\phie'$ preserves the same sign on $(-1,1)$, $\forall\epsilon>0$.
Without loss of generality, we may assume $\phie'(x)>0$ for $x\in(-1,1)$. Then by (\ref{2-eqn2}), one may get
$\phie(x)\geq\phie(-1)\geq\phi_0^-\geq\min\{\phi_0^+,\phi_0^-\}$.
Along with (\ref{id:2-w}), we obtain $0=\lim_{\epsilon\downarrow0}\phie(0)\geq\min\{\phi_0^+,\phi_0^-\}$,
which contradicts to the assumption $\min\{\phi_0^+,\phi_0^-\}>0$.
Consequently, there exist $\epsilon^*>0$ and $x_\epsilon^*\in(-1,1)$
such that $\phie'(x_\epsilon^*)=0$ as $0<\epsilon<\epsilon^*$.
As for the proof of Theorem \ref{2-thm2} (i), we may use (\ref{id:2-314}) and the fact that $\phie''\geq0$ on $[-1,1]$
to get $\phie'(x_1)<0<\phie'(x_2)$ for $x_1\in(-1,x_\epsilon^*)$ and $x_2\in(x_\epsilon^*,1)$.
Hence, $\phie$ attains the minimum value at an interior point $x_\epsilon^*\in(-1,1)$. This completes the proof of Claim 2.
\end{proof}
By Claim 2, we complete the proof of Theorem~\ref{2-thm:w} (ii).
Similarly, Theorem~\ref{2-thm:w} (iii) can also be proved.

We prove Theorem~\ref{2-thm:w} (iv) in two cases: \textbf{(I)}
$\phie''$ never changes sign on $[-1,1]$; \textbf{(II)}
$\phie''$ changes sign on $[-1,1]$. For the case (I), without loss
of generality, we may assume $\phie''\geq0$ on $[-1,1]$. Then
by~(\ref{2-id:ww1}), $\phie\geq0$ on $[-1,1]$ and the maximum
value of $\phie$ occurs at the boundary $x=\pm1$. Suppose
$\phie(1)=\max_{x\in[-1,1]}\phie(x)$. Then $\phie'(1)\geq0$.
Moreover, by the boundary condition (\ref{2-eqn2}), we get
$\phi_0^+=\phie(1)+\eta_\epsilon\phie'(1)\geq0$, which gives
$\phi_0^-\leq0$ due to $\min\{\phi_0^+,\phi_0^-\}\leq0$.
Consequently, $\eta_\epsilon\phie'(-1)=\phie(-1)-\phi_0^-\geq0$.
Hence by the assumption of $\phie''\geq0$ on $[-1,1]$, we have
$\phie'(x)\geq\phie'(-1)\geq0$ for $x\in[-1,1]$, i.e., $\phie$ is
monotone increasing on $[-1,1]$. Along with the boundary condition
(\ref{2-eqn2}), we have
$\phi_0^-\leq\phie(-1)\leq\phie(x)\leq\phie(1)\leq\phi_0^+$ for
$x\in [-1,1]$. Similarly, if
$\phie(-1)=\max_{x\in[-1,1]}\phie(x)$, then we obtain
$\phie'\leq0$ on $[-1,1]$ and
$\phi_0^+\leq\phie(x)\leq\phi_0^-$. This proves~(iii). For the
case (II) there exists $\hat{x}_\epsilon\in (-1,1)$ such that
$\phie''(\hat{x}_\epsilon)=0$. Then we may use (\ref{id:2-314})
and the same argument as in Theorem~\ref{2-thm2}~(i) to get Theorem~\ref{2-thm:w}~(iv).

It remains to prove Theorem~\ref{2-thm:w}(v). Using
(\ref{2-id:w2}), Theorem~\ref{2-thm:w}(i) and the similar argument
of Lemma~\ref{2-lem3}, we may obtain
$\lim_{\epsilon\downarrow0}C_\epsilon'=-\frac{1}{2}f(0)$ which
implies \begin{equation}\label{20110301}
\lim_{\epsilon\downarrow0}\left[(\phi_0^+-\phie(1))^2-\gamma^2(f(\phie(1))-f(0))\right]=0\,,
\end{equation}
by setting $x=1$ in (\ref{2-id:w2}) and using boundary
condition~(\ref{2-eqn2}) with
$\lim_{\epsilon\downarrow0}\frac{\eta_\epsilon}{\epsilon}=\gamma$.
Therefore, as for the proof of Theorem~\ref{2-thm5}~(iv), we may
use Theorem~\ref{2-thm:w} (i)-(iii) and~(\ref{20110301}) to get
Theorem~\ref{2-thm:w}(v) and complete the proof of
Theorem~\ref{2-thm:w}.\\

\begin{rem}\label{2-rk5}
If $\sum_{k=1}^{N_1}a_k\alpha_k\neq\sum_{l=1}^{N_2}b_l\beta_l$, we have
$\lim_{\epsilon\downarrow0}w(x)=r$ for all $x\in(-1,1)$, where $r$ is uniquely
determined by $f'(r)=0$. The proof is similar to the proof of Theorem 4.2 of \cite{LHLL}.\\
\end{rem}

\section{Non-electroneutral cases: Proof of Theorem~\ref{NE-thm}}\label{NE-sec}

\ \ \ \ In this section, we assume $0<\alpha<\beta$ and $\phi_0^-=\phi_0^+$.
To prove Theorem~\ref{NE-thm}, we need the following properties, which can be obtained from \cite{LHLL}.
\begin{itemize}
\item[\textbf{(P1)}] \textbf{Gradient estimates of $\phie$} (cf. Theorem 3.1, \cite{LHLL}):
 The unique solution $\phie$ is even and satisfies $\phie''\leq0$ on $[-1,1]$, and $\phie'(x_1)\geq0\geq\phie'(x_2)$
 for $x_1\in[-1,0)$ and $x_2\in(0,1]$.
 Moreover, $\phie'$ satisfies
\begin{equation}\label{2013-1013-4}
-\phie'(-1)=\phie'(1)=\frac{\alpha-\beta}{2\epsilon^2}<0,
\end{equation}
and
\begin{equation}
|\phie'(x)|\leq\frac{\beta-\alpha}{\epsilon^2}\left(e^{-\frac{\sqrt{2\alpha}(1+x)}{2\epsilon}}+e^{-\frac{\sqrt{2\alpha}(1-x)}{2\epsilon}}\right),
\quad\forall\, x\in(-1,1).\label{2013-1013-10}
\end{equation}
\item[\textbf{(P2)}] \textbf{Interior asymptotic behavior of $\phie$} (cf. Theorem~1.5, \cite{LHLL}):
For any compact subset $K$ of $(-1,1)$, there holds
\begin{equation}\label{2013-1013-9}
\sup_{0<\epsilon<1}\left|\phie(x)-\phie(\pm1)-\log\frac{1}{\epsilon^2}\right|<\infty,\quad\forall\,x\in K.
\end{equation}

\item[\textbf{(P3)}] \textbf{Estimates of $n_\epsilon$ and $p_\epsilon$}:  In \cite{LHLL}, we have established the following estimates (see (3.9), (3.15) and (3.37) of \cite{LHLL}):
\begin{eqnarray}
4\leq&&\int_{-1}^1e^{\phie(y)}dy\int_{-1}^1e^{-\phie(y)}dy\leq\frac{4\beta}{\alpha},\label{0909-1}\\
\frac{\alpha e^{\phie(0)}}{\int_{-1}^1e^{\phie(y)}dy}+&&\frac{\beta e^{-\phie(0)}}{\int_{-1}^1e^{-\phie(y)}dy}
+\frac{\epsilon^2}{4}\int_{-1}^1\phie'^2(y)dy=\frac{\alpha+\beta}{2}\label{0909-2}
\end{eqnarray}
and
\begin{equation}
\frac{(\alpha-\beta)^2}{8\epsilon^2}\leq\frac{\beta e^{-\phie(1)}}{\int_{-1}^1e^{-\phie(y)}dy}\leq\frac{(\alpha-\beta)^2}{8\epsilon^2}+\frac{\alpha+\beta}{2}.\label{0909-3}
\end{equation}
Using (\ref{2013-1013-3})
and the fact that $\phie(1)=\phie(-1)$ and $\phie'(0)=0$ (by (P1)), we can transform (\ref{0909-1})-(\ref{0909-3}) into
\begin{eqnarray}
\frac{\alpha^2}{4}\leq n_\epsilon(x)&&p_\epsilon(x)\leq\frac{\alpha\beta}{4},\quad\forall\,x\in[-1,1],\label{2013-1013-6}\\
n_\epsilon(0)+p_\epsilon(0&&)+\frac{\epsilon^2}{4}\int_{-1}^1\phie'^2(y)dy=\frac{\alpha+\beta}{2},\label{2013-1013-7}\\
\frac{(\alpha-\beta)^2}{8\epsilon^2}\leq p_\epsilon(1&&)=p_\epsilon(-1)\leq\frac{(\alpha-\beta)^2}{8\epsilon^2}+\frac{\alpha+\beta}{2},\label{2013-1013-8}
\end{eqnarray}
respectively.
\end{itemize}

Having (P1)-(P3) at hand, we are now in a position to prove Theorem~\ref{NE-thm}.\\

\textbf{Proof of Theorem~\ref{NE-thm}.}
Let $I_{\epsilon^\kappa}=[-1+\epsilon^\kappa,1-\epsilon^\kappa]$, where $0<\epsilon,\,\kappa<1$. For any $y\in I_{\epsilon^\kappa}$, we may use (\ref{2013-1013-10}) to get
\begin{equation}\label{2014-0910-11}
|\phie(y)-\phie(0)|\leq\frac{\beta-\alpha}{\epsilon^2}\left|\int_0^y\left(e^{-\frac{\sqrt{2\alpha}(1+x)}{2\epsilon}}+e^{-\frac{\sqrt{2\alpha}(1-x)}{2\epsilon}}\right)dx\right|\leq\frac{2\sqrt{2}(\beta-\alpha)}{\sqrt{\alpha}\epsilon}
e^{-\frac{\sqrt{2\alpha}}{\epsilon^{1-\kappa}}}.
\end{equation}
As a consequence, we have
\begin{equation}\label{2013-1013-11}
|\phie(x)-\phie(y)|\leq|\phie(x)-\phie(0)|+|\phie(0)-\phie(y)|\leq\frac{2\sqrt{2}(\beta-\alpha)}{\sqrt{\alpha}\epsilon}
e^{-\frac{\sqrt{2\alpha}}{\epsilon^{1-\kappa}}},
\end{equation}
for $x,y\in I_{\epsilon^\kappa}$.
Note that $\lim_{\epsilon\downarrow0}\frac{2\sqrt{2}(\beta-\alpha)}{\sqrt{\alpha}\epsilon}
e^{-\frac{\sqrt{2\alpha}}{\epsilon^{1-\kappa}}}=0$ for $0<\kappa<1$. Thus (\ref{2013-1013-11}) gives
\begin{equation}\label{2014-0829-6}
\lim_{\epsilon\downarrow0}\sup_{x,y\in I_{\epsilon^\kappa}}|\phie(x)-\phie(y)|=0.
\end{equation}

For $0<\epsilon<1$, we may set $x=0$ in (\ref{2013-1013-9}) and combine the result with (\ref{2014-0910-11}) to get
\begin{equation}\label{2013-1013-12}
\sup_{0<\epsilon<1}\left|\phi(y)-\phi(\pm1)-\log\frac{1}{\epsilon^2}\right|<\infty,\quad\forall\, y\in\, I_{\epsilon^\kappa}.
\end{equation}
To prove (\ref{2014-0810-1pm})-(\ref{thm2014-0829}), we need the following claim:\\

\textbf{Claim 3.}
\begin{itemize}
\item[(i)] At the boundary $x=\pm1$, we have
\begin{equation}\label{2014-0829-2}
\lim_{\epsilon\downarrow0}\frac{n_\epsilon(\pm1)}{\epsilon^{2-\tau}}=0\quad \mathrm{and}\quad\lim_{\epsilon\downarrow0}\epsilon^2p_\epsilon(\pm1)=\frac{(\alpha-\beta)^2}{8},
\end{equation}
for any $\tau>0$.

\item[(ii)] Assume $0<\epsilon<1$. Then there exists $\lambda_\epsilon(\kappa)>0$ such that $\lim_{\epsilon\downarrow0}\lambda_\epsilon(\kappa)=0$
and
\begin{equation}\label{2014-0830-1}
 \frac{\alpha}{2}-\lambda_\epsilon(\kappa)\leq n_\epsilon(x)\leq\frac{\sqrt{\alpha\beta}}{2},\quad \frac{\alpha}{2}\leq p_\epsilon(x)\leq\frac{\beta}{2}+\lambda_\epsilon(\kappa),
\end{equation}
for $x\in\,I_{\epsilon^\kappa}$. Moreover,
\begin{equation}\label{2013-1014-4}
\lim_{\epsilon\downarrow0}\sup_{x\in I_{\epsilon^\kappa}}|n_\epsilon(x)-n_\epsilon(0)|
=\lim_{\epsilon\downarrow0}\sup_{x\in I_{\epsilon^\kappa}}|p_\epsilon(x)-p_\epsilon(0)|=0.
\end{equation}
\end{itemize}
\begin{proof}
(\ref{2013-1013-6}) and (\ref{2013-1013-8}) give
$\frac{2\alpha^2\epsilon^2}{(\alpha-\beta)^2+4(\alpha+\beta)\epsilon^2}\leq n_\epsilon(-1)=n_\epsilon(1)\leq\frac{2\alpha\beta\epsilon^2}{(\alpha-\beta)^2}$.
This shows $\lim_{\epsilon\downarrow0}\frac{n_\epsilon(\pm1)}{\epsilon^{2-\tau}}=0$ for any $\tau>0$.
Along with (\ref{2013-1013-8}), we prove (\ref{2014-0829-2}).

By (P1) and (\ref{2-eqn2014}), we have
\begin{equation}\label{phie-0829}
\phie(0)=\max_{x\in[-1,1]}\phie(x)
\end{equation}
and
\begin{equation}\label{2013-1013-5}
n_\epsilon(x)-p_\epsilon(x)=\epsilon^2\phie''(x)\leq0,\quad\forall\,x\in[-1,1].
\end{equation}
Along with (\ref{2013-1013-6}), we obtain
\begin{equation}\label{2014-0829-3}
p_\epsilon(x)\geq\frac{\alpha}{2}\quad and\quad n_\epsilon(x)\leq\frac{\sqrt{\alpha\beta}}{2},\quad\forall x\in[-1,1].
\end{equation}
By (\ref{2013-1013-3}), (\ref{2013-1013-7}) and (\ref{phie-0829}), one may check that
\begin{equation}\label{2013-1014-3}
0\leq n_\epsilon(0)-n_\epsilon(x)
=n_\epsilon(0)\left(1-e^{\phie(x)-\phie(0)}\right)
\leq\frac{\alpha+\beta}{2}\left(1-e^{\phie(x)-\phie(0)}\right),
\end{equation}
and
\begin{equation}\label{2013-1014-2}
0\leq p_\epsilon(x)-p_\epsilon(0)
=p_\epsilon(0)\left(e^{-\phie(x)+\phie(0)}-1\right)
\leq\frac{\alpha+\beta}{2}\left(e^{-\phie(x)+\phie(0)}-1\right).
\end{equation}
Consequently, by (\ref{2014-0829-6}), (\ref{2013-1014-3}) and (\ref{2013-1014-2}),
we get (\ref{2013-1014-4}).

It remains to prove (\ref{2014-0830-1}). Let
\begin{equation}\label{2014-123}
\ds\lambda_\epsilon(\kappa)=\max\left\{\sup_{x\in I_{\epsilon^\kappa}}|n_\epsilon(x)-n_\epsilon(0)|,
\sup_{x\in\,I_{\epsilon^\kappa}}|p_\epsilon(x)-p_\epsilon(0)|\right\}>0.
\end{equation}
By (\ref{2013-1014-4}), we have $\lim_{\epsilon\downarrow0}\lambda_\epsilon(\kappa)=0$.
Using (\ref{phie-0829}), one may find
\begin{equation}
n_\epsilon(0)\equiv\frac{\alpha e^{\phie(0)}}{\int_{-1}^1e^{\phie(y)} dy}=\frac{\alpha}{\int_{-1}^1e^{\phie(y)-\phie(0)} dy}\geq\frac{\alpha}{2},\label{2014-0830-2}
\end{equation}
and
\begin{equation}
p_\epsilon(0)\equiv\frac{\beta e^{-\phie(0)}}{\int_{-1}^1e^{-\phie(y)} dy}=\frac{\beta}{\int_{-1}^1e^{-\phie(y)+\phie(0)} dy}\leq\frac{\beta}{2}.\label{2014-0830-3}
\end{equation}
Hence, (\ref{2014-0829-3}) and (\ref{2014-0830-2}) immediately give $\frac{\sqrt{\alpha\beta}}{2}\geq n_\epsilon(x)\geq
n_\epsilon(0)-\lambda_\epsilon(\kappa)\geq\frac{\alpha}{2}-\lambda_\epsilon(\kappa)$, for $x\in\, I_{\epsilon^{\kappa}}$.
On the other hand, by (\ref{2014-0829-3}) and (\ref{2014-0830-3}) we obtain $\frac{\alpha}{2}\leq p_\epsilon(x)\leq p_\epsilon(0)+
\lambda_\epsilon(\kappa)\leq\frac{\beta}{2}+\lambda_\epsilon(\kappa)$, for $x\in\, I_{\epsilon^{\kappa}}$. Therefore, we get (\ref{2014-0830-1}) and
 complete the proof of Claim~3.
\end{proof}

(\ref{2014-0810-1pm}) immediately follows from (\ref{2013-1013-5}) and (\ref{2014-0830-1}), and (\ref{2014-0910-2}) follows from (\ref{2014-0829-2}).
To prove (\ref{thm2014-0829}), we rewrite $n_\epsilon(0)=\frac{\alpha e^{\phie(0)}}{\int_{-1}^1e^{\phie(y)} dy}$ as
\begin{equation}\label{2013-1013-13}
n_\epsilon(0)=\frac{\alpha}
{\left(\int_{-1}^{-1+\epsilon^\kappa}+\int_{-1+\epsilon^\kappa}^{1-\epsilon^\kappa}+\int_{1-\epsilon^\kappa}^1\right)
e^{\phie(y)-\phie(0)} dy}.
\end{equation}
By (\ref{phie-0829}), we have
\begin{equation}\label{2013-1013-14}
0\leq\left(\int_{-1}^{-1+\epsilon^\kappa}+\int_{1-\epsilon^\kappa}^1\right)e^{\phie(y)-\phie(0)}dy\leq2\epsilon^\kappa.
\end{equation}
On the other hand, by (\ref{2014-0829-6}) we get
\begin{equation}\label{2013-1014-1}
\lim_{\epsilon\downarrow0}\int_{-1+\epsilon^\kappa}^{1-\epsilon^\kappa}e^{\phie(y)-\phie(0)}dy=2.
\end{equation}
Combining (\ref{2013-1013-13})-(\ref{2013-1014-1}), we conclude that
\begin{equation}\label{2013-1013-15}
\lim_{\epsilon\downarrow0}n_\epsilon(0)=\frac{\alpha}{2}.
\end{equation}

To deal with the limit of value $p_\epsilon(0)$ as $\epsilon$ tends to zero,
 we need the following estimate:
\begin{eqnarray}\label{2014-246}
\left|p_\epsilon(0)-\frac{\alpha}{2}\right|&\leq|n_\epsilon(0)-p_\epsilon(0)|+\left|n_\epsilon(0)-\frac{\alpha}{2}\right|\nonumber\\[-0.7em]
&\\[-0.7em]
&\leq
 |n_\epsilon(x)-p_\epsilon(x)|+\left|n_\epsilon(0)-\frac{\alpha}{2}\right|
+2\lambda_\epsilon(\kappa),\quad\forall x\in I_{\epsilon^\kappa}. \nonumber
\end{eqnarray}
Here we have used (\ref{2013-1014-4}) and (\ref{2014-123}) to get the second line of (\ref{2014-246}).
On the other hand, by integrating (\ref{2-eqn2014}) over $I_{\epsilon^\kappa}$ and using (\ref{2013-1013-3}) and (\ref{2013-1013-10}),
we obtain
\begin{equation}\label{2013-1013-17}
0\leq\int_{-1+\epsilon^\kappa}^{1-\epsilon^\kappa}(p_\epsilon(x)-n_\epsilon(x))dx=\epsilon^2\left(\phie'(-1+\epsilon^\kappa)-\phie'(1-\epsilon^\kappa)\right)
\leq 4(\beta-\alpha)e^{-\frac{\sqrt{2\alpha}}{2\epsilon^{1-\kappa}}}.
\end{equation}
Note that $0<\kappa<1$. As a consequence, by (\ref{2013-1013-15})-(\ref{2013-1013-17}) we find
\begin{equation}\label{2013-1013-18}
\lim_{\epsilon\downarrow0}p_\epsilon(0)=\frac{\alpha}{2}.
\end{equation}
Then (\ref{thm2014-0829}) follows from (\ref{2013-1014-4}), (\ref{2013-1013-15}) and (\ref{2013-1013-18}).

By (\ref{2014-0829-3}), we immediately get (\ref{thm1-6}).
Now we shall prove (\ref{thm1-7}).
Let $g(x)\in C^1([-1,1])$.
Multiplying (\ref{2-eqn2014}) by $g(x)$ and
integrating the result over $(-1,1)$, we have
\begin{eqnarray}\label{2013-1014-5}
\int_{-1}^1(n_\epsilon(x)-p_\epsilon(x))g(x)dx=&\epsilon^2\int_{-1}^1\phie''(x)g(x)dx\nonumber\\[-0.7em]
&\\[-0.7em]
=&\frac{\alpha-\beta}{2}\left(g(-1)+g(1)\right)-\epsilon^2\int_{-1}^1\phie'(x)g'(x)dx.\nonumber
\end{eqnarray}
Here we have used the intergration by parts and (\ref{2013-1013-4}) to get (\ref{2013-1014-5}).
On the other hand, by using (\ref{2013-1013-10}), one may check that
\begin{eqnarray}\label{2013-1014-6}
\left|\epsilon^2\int_{-1}^1\phie'(x)g'(x)dx\right|\leq&&
(\beta-\alpha)\max_{x\in[-1,1]}|g(x)|\int_{-1}^1\left(e^{-\frac{\sqrt{2\alpha}(1+x)}{2\epsilon}}+e^{-\frac{\sqrt{2\alpha}(1-x)}{2\epsilon}}\right)dx\nonumber\\[-0.7em]
\\[-0.7em]
\leq&&2(\beta-\alpha)\sqrt{\frac{2}{\alpha}}\left(\max_{x\in[-1,1]}|g(x)|\right)\epsilon.\nonumber
\end{eqnarray}
By  (\ref{thm1-6}), (\ref{2013-1013-5}), (\ref{2014-0829-3}), (\ref{2013-1013-17}), (\ref{2013-1014-5}) and (\ref{2013-1014-6}), we have
\begin{eqnarray}
&&\left|\left(\int_{-1}^{-1+\epsilon^\kappa}+\int^{1}_{1-\epsilon^\kappa}\right)p_\epsilon(x)g(x)dx-\frac{\beta-\alpha}{2}\left(g(-1)+g(1)\right)\right|\nonumber\\
=&&\left|\left(\int_{-1}^{-1+\epsilon^\kappa}+\int^{1}_{1-\epsilon^\kappa}\right)n_\epsilon(x)g(x)dx
+\int_{-1+\epsilon^\kappa}^{1-\epsilon^\kappa}(n_\epsilon(x)-p_\epsilon(x))g(x)dx+\epsilon^2\int_{-1}^1\phie'(x)g'(x)dx\right|\nonumber\\
\leq&&\max_{x\in[-1,1]}|g(x)|\times\left[\sqrt{\alpha\beta}\epsilon^\kappa+4(\beta-\alpha)e^{-\frac{\sqrt{2\alpha}}{2\epsilon^{1-\kappa}}}
+2(\beta-\alpha)\sqrt{\frac{2}{\alpha}}\epsilon\right].\nonumber
\end{eqnarray}
Note that $0<\kappa<1$. Consequently,
\begin{equation}\label{2013-1014-7}
\lim_{\epsilon\downarrow0}\left(\int_{-1}^{-1+\epsilon^\kappa}+\int^{1}_{1-\epsilon^\kappa}\right)p_\epsilon(x)g(x)dx=\frac{\beta-\alpha}{2}\left(g(-1)+g(1)\right).
\end{equation}
In particular, let $g\in C^1([-1,1])$ satisfy $g(x)=1$ for $x\in[-1,0]$, $g(x)\in[0,1]$ for $x\in[0,1/2]$, and $g(x)=0$ for $x\in[1/2,1]$.
Then (\ref{2013-1014-7}) gives $\lim_{\epsilon\downarrow0}\int_{-1}^{-1+\epsilon^\kappa}p_\epsilon(x)dx=\frac{\beta-\alpha}{2}$.
Similarly, we have $\lim_{\epsilon\downarrow0}\int^{1}_{1-\epsilon^\kappa}p_\epsilon(x)dx=\frac{\beta-\alpha}{2}$.
Therefore, we get (\ref{thm1-7}) and complete the proof of Theorem~\ref{NE-thm}(i).

It remains to prove (\ref{2013-1015-1}). By (\ref{2013-1013-3}), we have
\begin{equation}\label{2014-0910-9}
\phie(x)-\phie(1)-\log\frac{1}{\epsilon^2}=\log\frac{\epsilon^2p_\epsilon(1)}{p_\epsilon(x)}.
\end{equation}
Note that for any compact subset $K$ of $(-1,1)$, we have
$K\subset I_{\epsilon^\kappa}$ as $0<\epsilon\ll1$ is sufficiently small.
Hence, by (\ref{thm2014-0829}), (\ref{2014-0829-2}) and (\ref{2014-0910-9}), we conclude that
$\lim_{\epsilon\downarrow0}\left(\phie(x)-\phie(1)-\log\frac{1}{\epsilon^2}\right)=\log\frac{\lim_{\epsilon\downarrow0}\epsilon^2p_\epsilon(1)}{\lim_{\epsilon\downarrow0}p_\epsilon(x)}=\log\frac{(\alpha-\beta)^2}{4\alpha}$
uniformly in $K$. Therefore, we get (\ref{2013-1015-1}) and complete the proof of Theorem~\ref{NE-thm}.\\

\section{Numerical experiments}\label{2-mpnpsec:5}
\medskip
\noindent

In this section, we do numerical
computations to compare solutions
of the CCPB and PB equations. All numerical results are obtained
using the convex iteration method \cite{LHLL, R1, RLZ, RLW} and
the finite element methods with piecewise linear space
which is used to solve the linearized equations. The
computational domain and the mesh size $h$ are fixed with is $[ -1, 1 ]$, $h =
2^{-11}$, respectively, throughout the numerical experiments. The
values of $\epsilon$ are set by $\epsilon = 2^{-j}, j = 1,3,5$, in order to observe
the tendency of the associated solutions $\phie$'s as $\epsilon$
goes to zero.

As for~\cite{LHLL}, the numerical
scheme can be extended to the CCPB equation (\ref{2-eqn1}) with
the boundary condition (\ref{2-eqn2}) and it can be presented as
follows:

\begin{eqnarray}\label{eqn:convex1}
   \epsilon^2 \phi_{m+\frac12}'' & = & \displaystyle \sum_{k=1}^{N_1} \frac{a_k \alpha_k}{\int_{-1}^1 e^{a_k \phi_m}\,dx} e^{a_k \phi_m} -\sum_{l=1}^{N_2} \frac{b_l \beta_l}{\int_{-1}^1 e^{-b_l \phi_m}\,dx} e^{-b_l \phi_m},\\
  \phi_{m+1} & = & s \phi_{m+\frac12} +( 1 -s ) \phi_m\,,
  \label{convex1.1}
  \end{eqnarray}
for $m = 1, 2, \cdots$, where $s$ is a positive constant
satisfying $0<s < 1$ with boundary conditions
\begin{eqnarray}\label{eqn:convex2}
 \phi_{m+\frac12}(-1) -\eta_{\epsilon} \phi_{m+\frac12}'(-1) = \phi_0^-, \quad \phi_{m+\frac12}(1) +\eta_{\epsilon} \phi_{m+\frac12}'(1) = \phi_0^+.
\end{eqnarray}
Let $\phi_{m+\frac12} = \phi_m +\delta_m$ with the correction term
$\delta_m$ which satisfies
\begin{eqnarray}\label{eqn:convex3}
 \delta_m(-1) -\eta_{\epsilon} \delta_m'(-1) = 0, \quad \delta_m(1) +\eta_{\epsilon} \delta_m'(1) =
 0
\end{eqnarray}
so that $\phi_{m+1} = \phi_m +s \delta_m = \phi_1 +s \sum_{i=1}^m
\delta_i$. If $\lim_{m \rightarrow \infty} |\delta_m| = 0$, then
the iterative scheme converges.

Define the residual function
$\mathcal{R}(\phi_m)$ as
\begin{equation}\label{eqn:convex4}
 \mathcal{R}(\phi_m) = \sum_{k=1}^{N_1} \frac{a_k \alpha_k}{\int_{-1}^1 e^{a_k \phi_m}\,dx} e^{a_k \phi_m}
 -\sum_{l=1}^{N_2} \frac{b_l \beta_l}{\int_{-1}^1 e^{-b_l \phi_m}\,dx} e^{-b_l \phi_m} -\epsilon^2
 \phi_m''\,.
\end{equation}
Then we obtain
\begin{eqnarray}
 \epsilon^2 \phi_{m+1}'' -\epsilon^2 \phi_m'' = s \epsilon^2 \delta_m'' = s \mathcal{R}(\phi_m).\label{eqn:convex5}
\end{eqnarray}
Integrating $\mathcal{R}(\phi_{m+1}) -\mathcal{R}(\phi_m)$, we may
use (\ref{eqn:convex4}) and (\ref{eqn:convex5}) to get
\begin{eqnarray}\label{eqn:convex6}
 \int^1_{-1} \mathcal{R}(\phi_{m+1})\,dx =  ( 1 -s ) \int^1_{-1} \mathcal{R}(\phi_m)\,dx .
\end{eqnarray}
In case of $s = 1$ in (\ref{eqn:convex6}), numerical scheme may
not converge and oscillate during the iteration procedure. When $0
< s < 1$, we have empirically observed that the value of $s$
should be compatible to $C \epsilon^2$ in order to let the
iteration converge. Moreover, the value of $C$ is chosen in the
interval $(0, 1)$ so that the convergence of the scheme can be
guaranteed. In the iteration procedure, the value $10^{-6}$ is
applied for stopping criterion with $|| \delta_m ||_{\infty} = ||
( \phi_{m+1} -\phi_m )/s ||_{\infty}$.

\begin{figure}[hbt]
 \centering{%
  \begin{tabular}{@{\hspace{-0.7pc}}c@{\hspace{-1.2pc}}c}
  \includegraphics[width=2.9in]{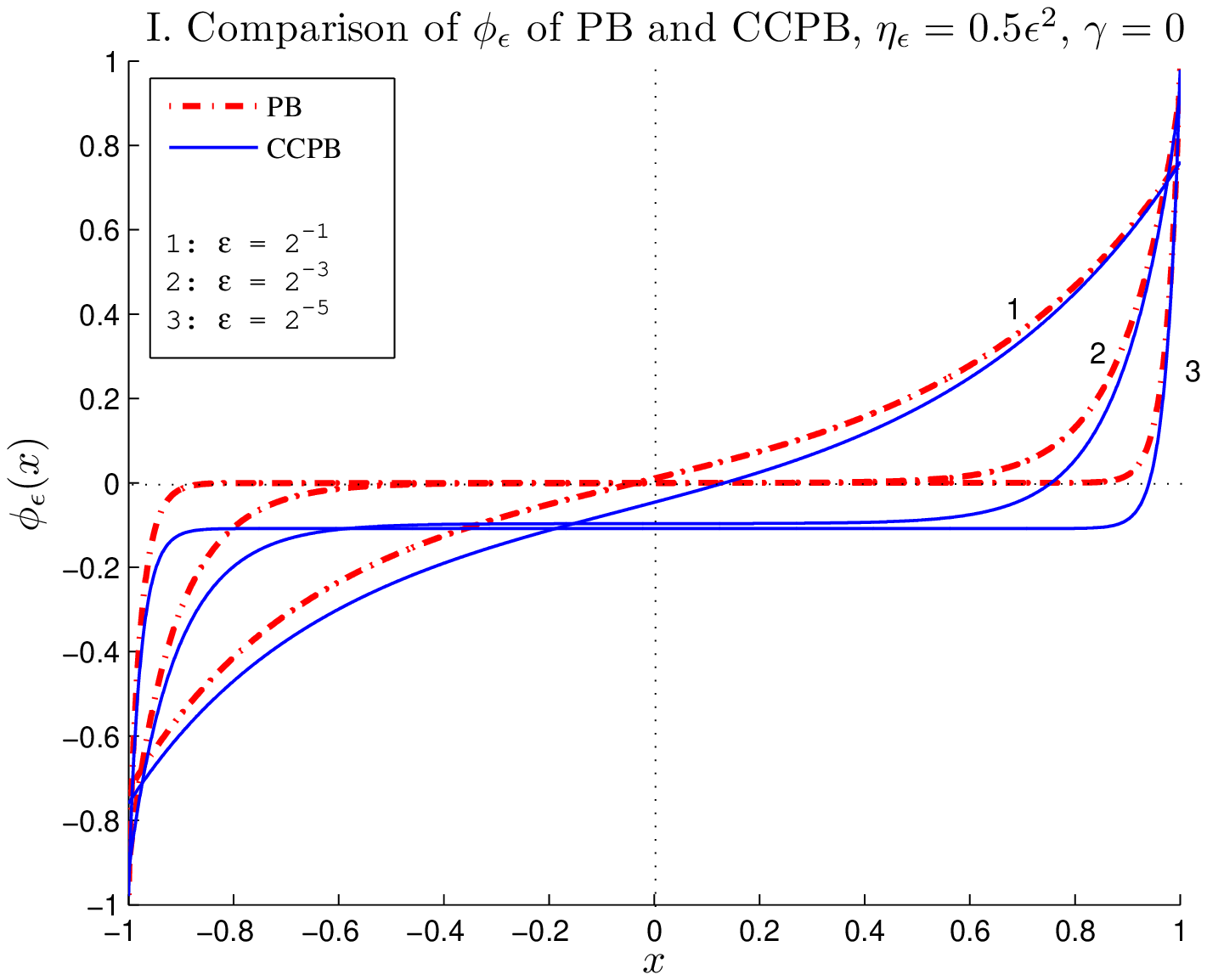} & \includegraphics[width=2.9in]{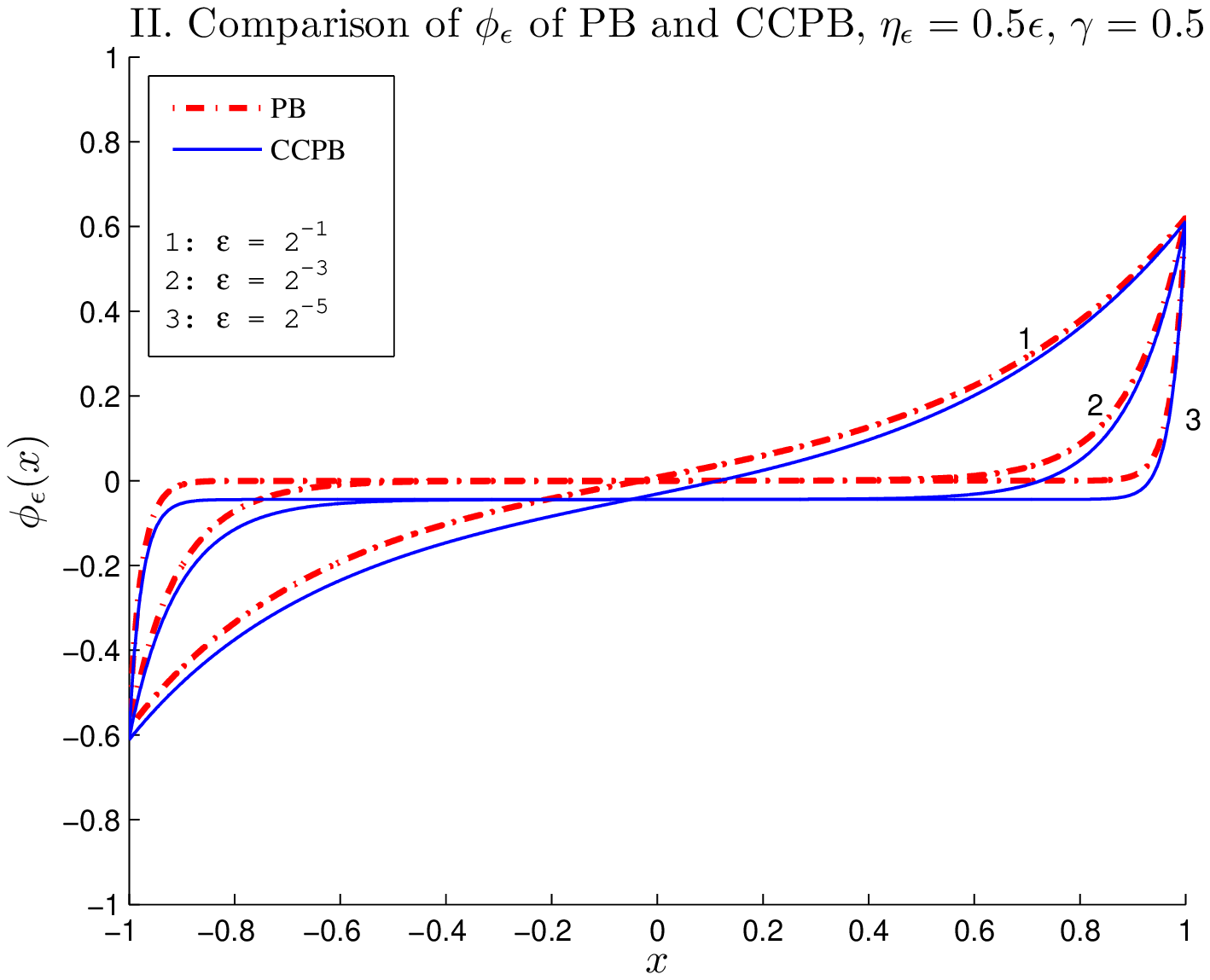}
  \end{tabular}}
 \caption{\small Comparison of $\phie$ of PB and CCPB equation with the electroneutral condition. $\phie$ of PB equation are in (red) dash-dotted curves and $\phie$ of CCPB equation are in (blue) solid curves. The label of curves in each picture depends on the dielectric constant $\epsilon = 2^{-1}, 2^{-3}, 2^{-5}$.
  I. $\eta_{\epsilon} = 0.5 \epsilon^2$ and $\gamma = 0$. II. $\eta_{\epsilon} = 0.5 \epsilon$ and $\gamma = 0.5$. In this computations, three ion species are used, one anion with valence $-1$, $\alpha_1=1.2$ and two cations with valences $1,2$
 , $\beta_1=\beta_2=0.4$. }\label{fig:F_1}
\end{figure}

For the PB equation
(\ref{2-eqnw}), we replace the denominators of the right hand
side of the equation (\ref{eqn:convex1}) by the value $2$. Then as
for the scheme of (\ref{eqn:convex1})-(\ref{eqn:convex2}), we
have a similar way to solve the PB equation~(\ref{2-eqnw}) with
the boundary condition (\ref{2-eqn2}), numerically. To compare
solutions of the PB and CCPB equations, we firstly set the
parameters as $N_1=1, N_2=2$, $a_1 = b_1 = 1, b_2 = 2$ and
$\alpha_1 = 1.2, \beta_1 = \beta_2 = 0.4$ so that the
electroneutral condition $a_1\alpha_1 = b_1\beta_1 + b_2\beta_2$
holds. The numerical computations also impose the boundary data as
$\phi_0^+ = -\phi_0^- = 1$ and the values $\eta_{\epsilon}$'s
for the boundary conditions (\ref{2-eqn2}) as $\eta_{\epsilon} =
0.5 \epsilon^2$ and $0.5 \epsilon$ which include the cases of
$\lim_{\epsilon \downarrow 0} \frac{\eta_{\epsilon}}{\epsilon} =0$
and $0.5$. The corresponding results are presented in
Figure~\ref{fig:F_1} and Table~\ref{table:D_c} consistent with
Theorem~\ref{2-thm5} and~\ref{2-thm:w}.

In Figure~\ref{fig:F_1}, one may
see the difference between the solutions of (\ref{2-eqn1}) and
(\ref{2-eqnw}) with the same boundary condition (\ref{2-eqn2}) and the valence $z_i=-1$ for the anion, $i=1$, and
$z_j=1,2$ for the cations, $j=1,2$, respectively.
The solution profiles of the PB equation (\ref{2-eqnw}) are
plotted as (red) dash-dotted curves and those of the CCPB
equation (\ref{2-eqn1}) are sketched as (blue) solid curves. Here
the index numbers, $1,2,3$ are associated with various values of
$\epsilon$'s, and a (black) dotted line is represented as the
axes for a reference.

Table~\ref{table:D_c} shows the
numerical results of $\phie(0)$ and $c$ for the CCPB and PB
equations where the value $c$ is defined in Theorem~\ref{2-thm5} can
be computed by Newton's method. One can easily see that for the PB
equation, the value $c$ is always equal to zero but for the CCPB
equation, the value $c$ may not be equal to zero. The ratio
$\beta_1/\beta_2$ may affect the value $c$ and $t$. As
$\beta_1/\beta_2$ varies, the numerical values of~$\phie(0)$,
$\phie(1)$, $c$ and $t$ are presented in Table~\ref{table:D_1} for
the case of $a_1=b_1=1$, $b_2=2$ and $\epsilon=2^{-5}$. Note that
the numerical values of $\phie(0)$ and $\phie(1)$ are quite close
to those of $c\in (c^{\ast}, 0)$ and $t$, respectively. This is
consistent with the results of Theorem~\ref{2-thm5}.
We remark that if $t$ is fixed and $\beta_1/\beta_2$ is decreasing, then the value $c$ is
decreasing.

\begin{table}[thb]
 \begin{center}\caption{The numerical results of $\phie(0)$ and its limit value $c$ of PB and CCPB equation in Figure \ref{fig:F_1}.}\label{table:D_c}
  \begin{tabular}{c|c||c|c|c|c}
   \hline
    & $\epsilon$ & $2^{-1}$ & $2^{-3}$ & $2^{-5}$ & c\\
   \hline
   \hline
   I & PB    &  0.0106 &  0.0000 &  0.0000 & 0\\
     & CCPB & -0.0459 & -0.0964 & -0.1081 & -0.1126\\
   \hline
   \hline
   II & PB    &  0.0079 &  0.0000 &  0.0000 & 0\\
      & CCPB & -0.0311 & -0.0442 & -0.0442 & -0.0441\\
   \hline
  \end{tabular}
 \end{center}
\end{table}

\begin{table}[hbt]
 \begin{center}\caption{The numerical results of $\phie(1)$, $\phie(0)$ of CCPB equation and their limit values $t$, $c$, $c_{\ast}$ in (\ref{ca1}) where
 $\alpha_1 = \beta_1 +2 \beta_2$. $\beta_1 \sim [\rm{Na^+}]$ is fixed to $1$.
 $\epsilon$ is fixed to $2^{-5}$.}\label{table:D_1}
  \begin{tabular}{c||c|c|c|c|c|c}
   \hline
   $\eta_{\epsilon}$ & $\beta_1/\beta_2$ & $\phie(1)$ & $t$ & $\phie(0)$ & $c$ & $c_{\ast}$ \\
   \hline
   \hline
                   & 1   & 1.0000 & 1.0000 & -0.1124 & -0.1126 & -0.1446 \\
   $0$             & 1/2 & 1.0000 & 1.0000 & -0.1265 & -0.1265 & -0.1446 \\
                   & 1/3 & 1.0000 & 1.0000 & -0.1320 & -0.1320 & -0.1446 \\
   \hline
                   & 1   & 0.9679 & 1.0000 & -0.1059 & -0.1126 & -0.1446 \\
   $0.5\epsilon^2$ & 1/2 & 0.9581 & 1.0000 & -0.1171 & -0.1265 & -0.1446 \\
                   & 1/3 & 0.9504 & 1.0000 & -0.1206 & -0.1320 & -0.1446 \\
   \hline
                   & 1   & 0.4962 & 0.4960 & -0.0299 & -0.0299 & -0.0394 \\
   $0.5\epsilon$   & 1/2 & 0.4278 & 0.4277 & -0.0255 & -0.0255 & -0.0296 \\
                   & 1/3 & 0.3853 & 0.3853 & -0.0218 & -0.0218 & -0.0242 \\
   \hline
  \end{tabular}
 \end{center}
\end{table}

\begin{figure}[hbt]
 \centering{%
  \begin{tabular}{@{\hspace{-0.7pc}}c@{\hspace{-1.2pc}}c}
   \includegraphics[width=2.9in]{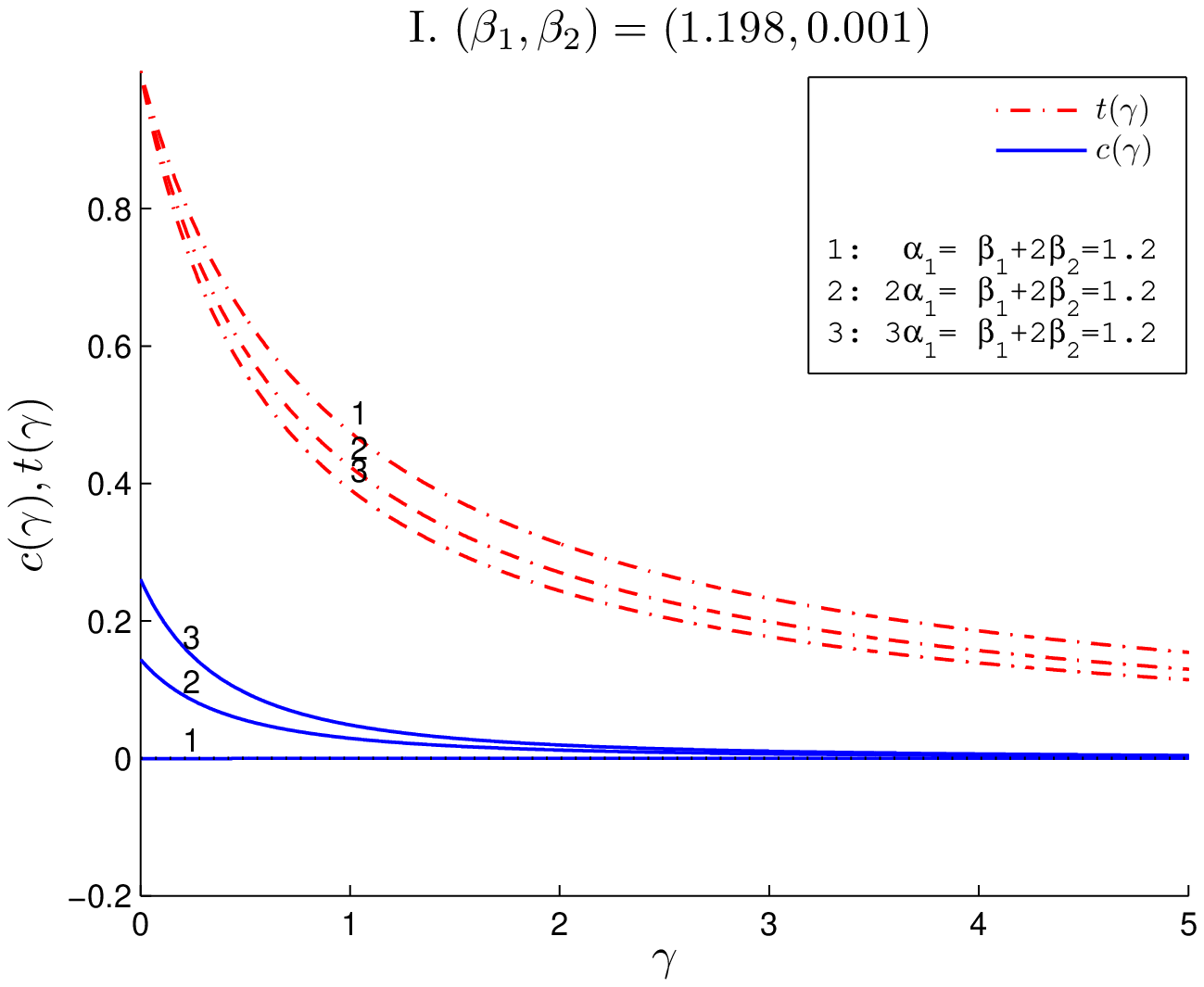} & \includegraphics[width=2.9in]{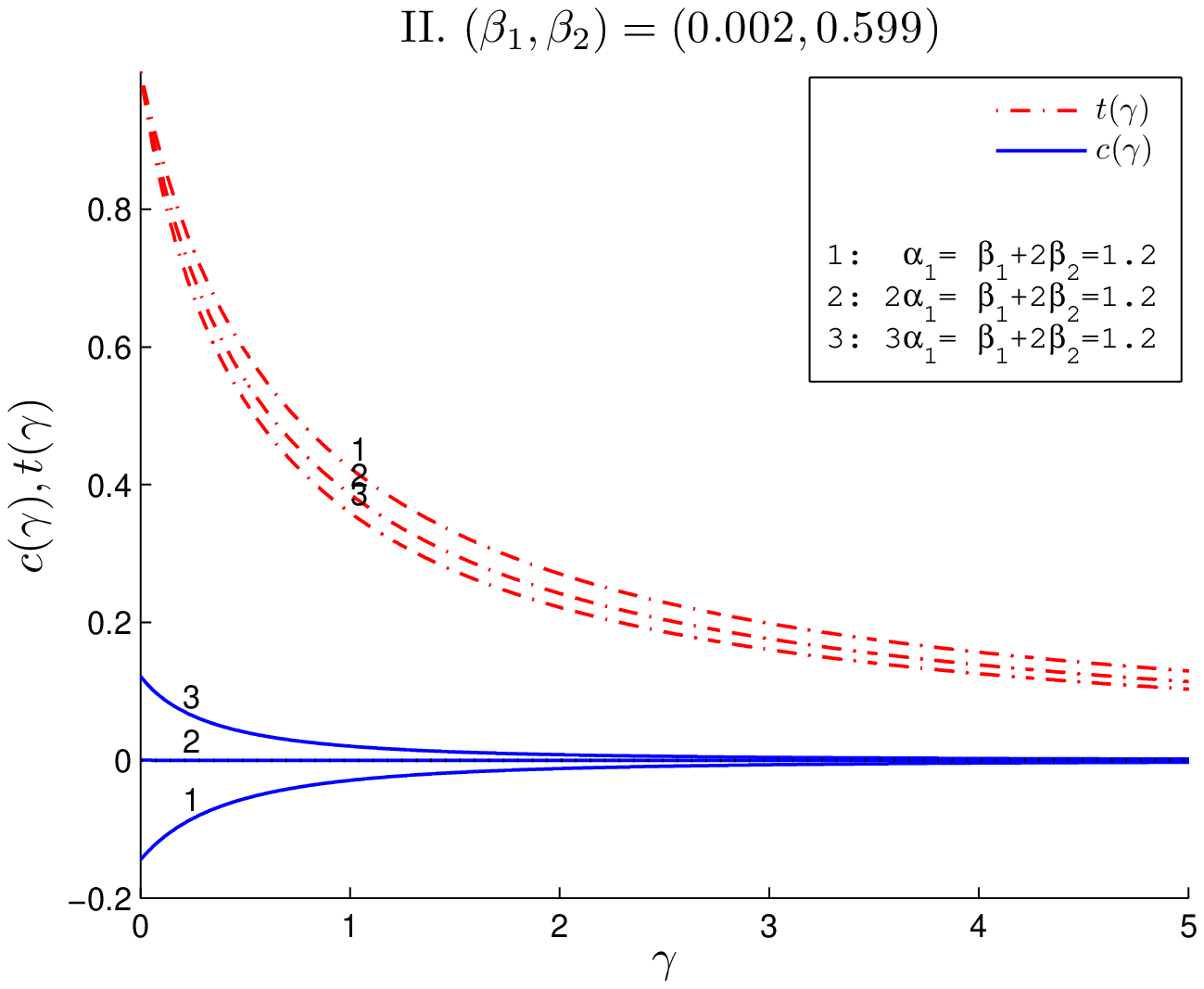}
  \end{tabular}}
 \caption{\small Comparison of $c(\gamma)$, $t(\gamma)$ with three species;
 one negative charge, two positive charges where $\alpha_1=1.2,0.6,0.4$ for $1$, $2$, $3$, respectively.
 I.   $(\beta_1,\beta_2) =(1.199,0.0005)$.
 II.  $(\beta_1,\beta_2) =(0.002,0.599)$.}\label{fig:ct_3_12_allm}
\end{figure}

\begin{figure}[hbt]
 \centering{%
  \begin{tabular}{@{\hspace{-0.7pc}}c@{\hspace{-1.2pc}}c}
   \includegraphics[width=2.9in]{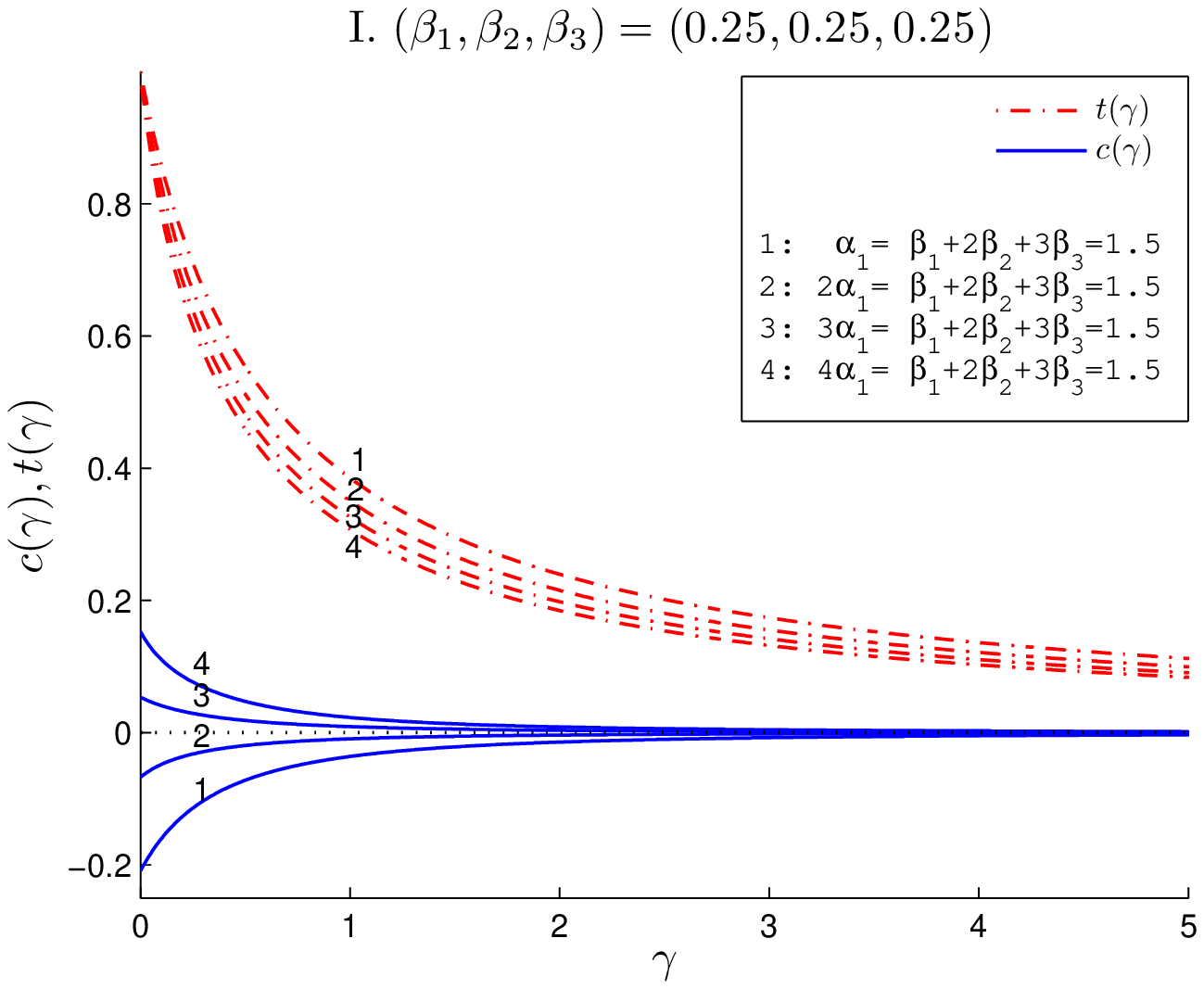} & \includegraphics[width=2.9in]{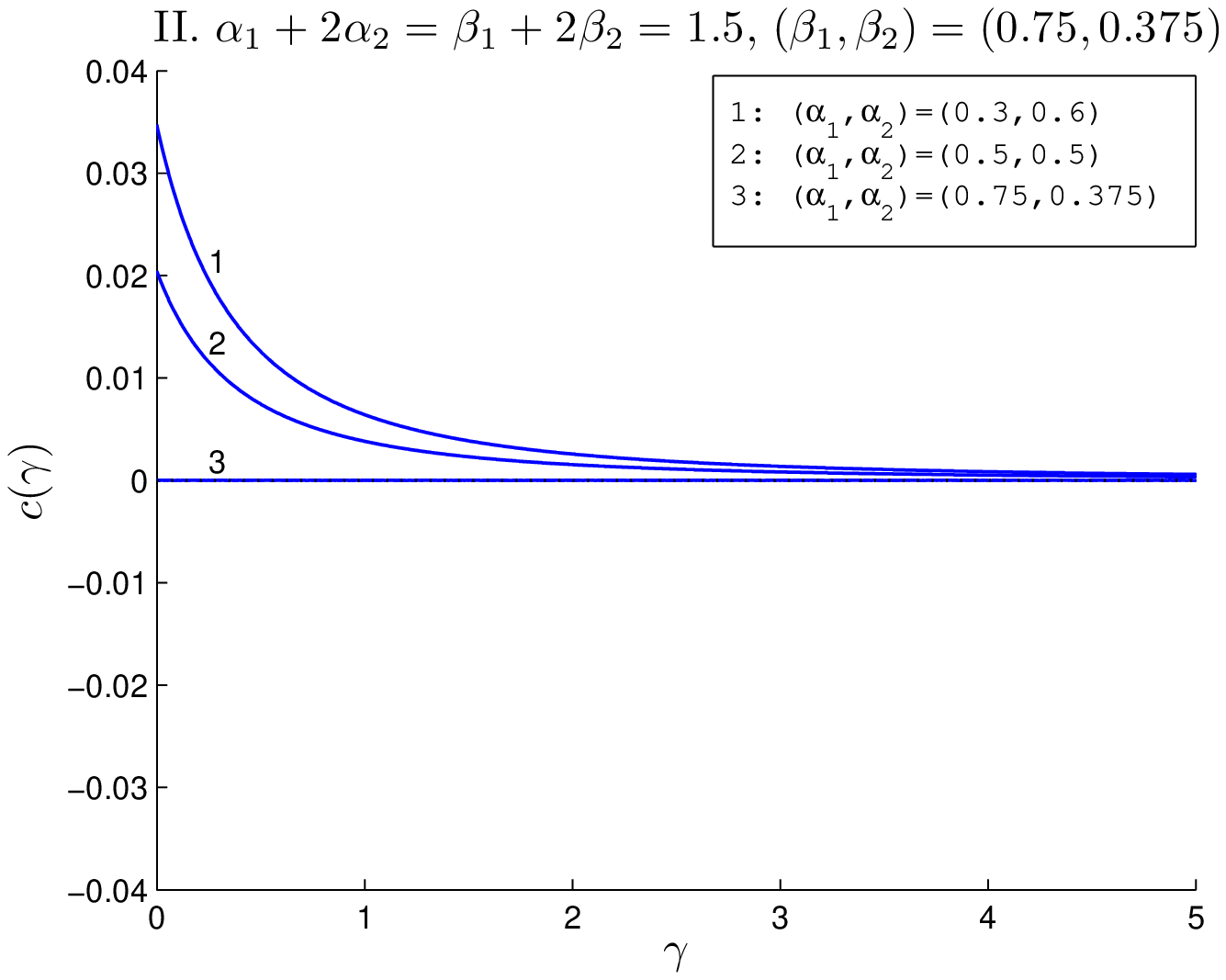}\\
   \includegraphics[width=2.9in]{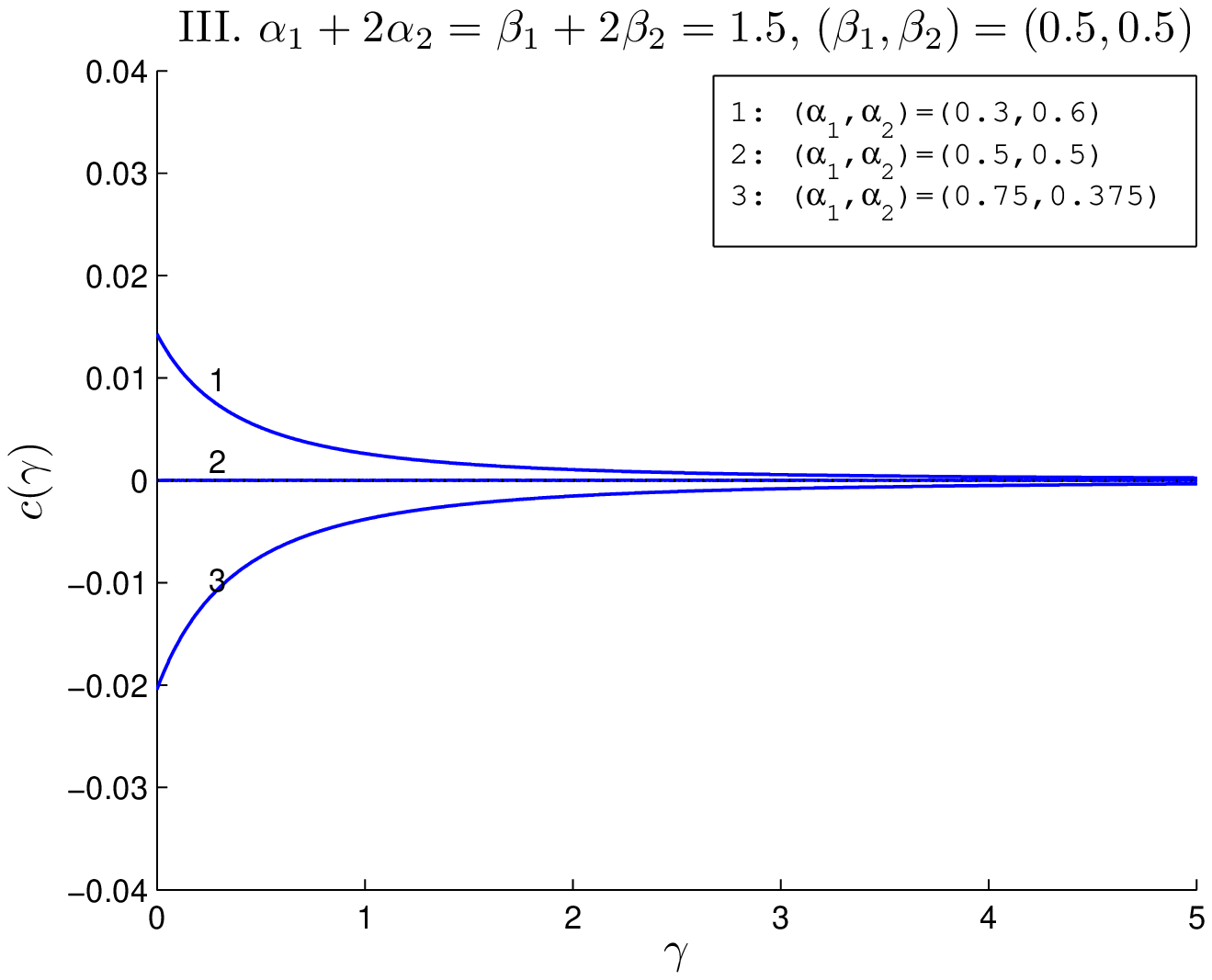} & \includegraphics[width=2.9in]{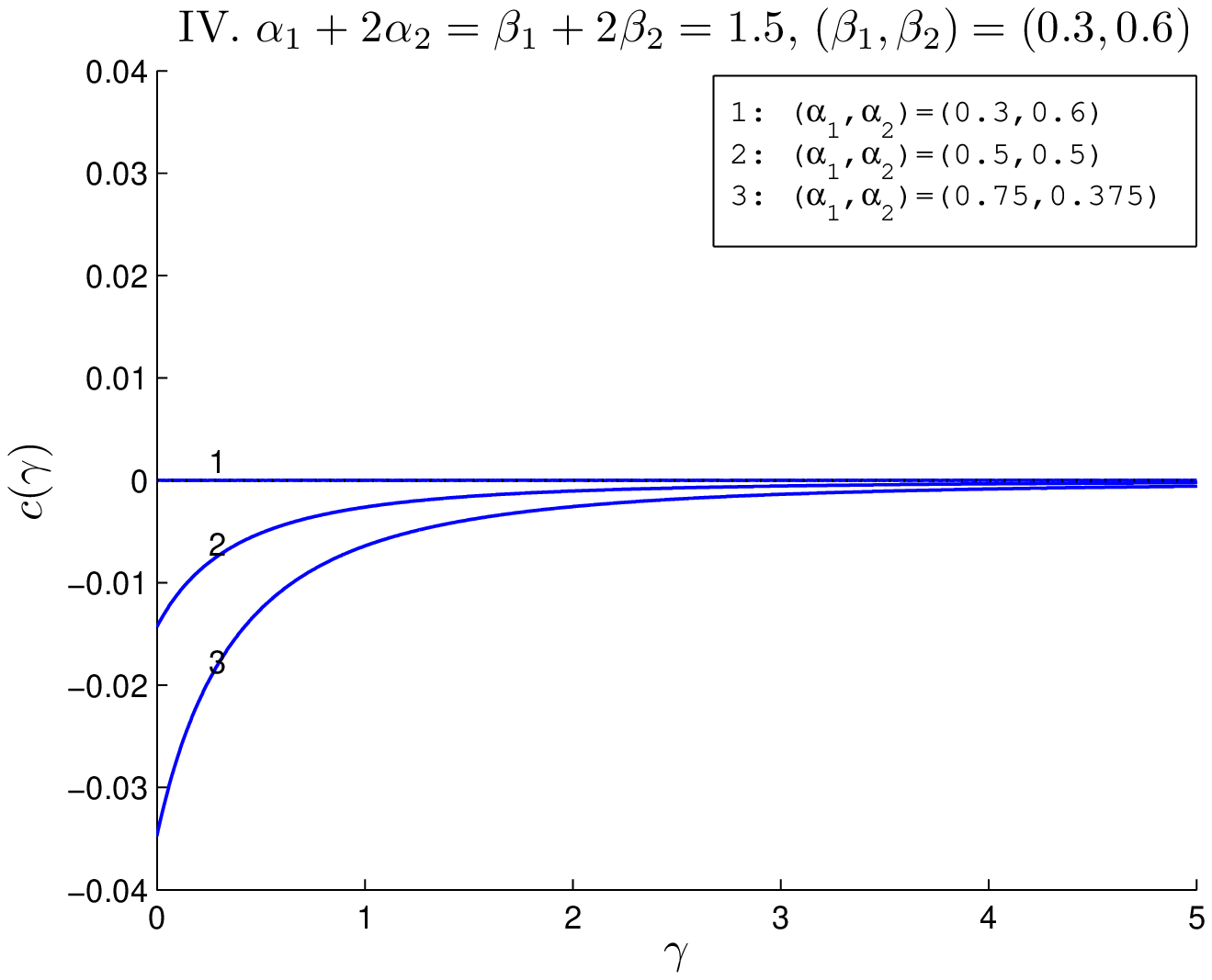}
  \end{tabular}}
 \caption{\small Comparison of $c(\gamma)$, $t(\gamma)$ with four ion species;
 one negative charges,
 three positive charges (I), and two negative charge, two positive charges (II, III, IV).
 I. $\alpha_1=1.5,0.75,0.5,0.375$ for $1$, $2$, $3$, $4$, respectively,
 and $(\beta_1,\beta_2,\beta_3)=(0.25,0.25,0.25)$.
 II. $(\beta_1,\beta_2)=(0.75,0.375)$.
 III. $(\beta_1,\beta_2)=(0.5,0.5)$.
 IV. $(\beta_1,\beta_2)=(0.3,0.6)$.
 ($1. (\alpha_1,\alpha_2)=(0.3,0.6)$,
  $2. (\alpha_1,\alpha_2)=(0.5,0.5)$,
  $3. (\alpha_1,\alpha_2)=(0.75,0.375)$ for II, III, IV.)}\label{fig:ct_4_all}
\end{figure}

From
Theorem~\ref{2-thm5}~(ii)--(iv), both $t$ and $t-c$ are decreasing
functions to $\gamma$. Surely, $c$ can be regarded as a function
to $\gamma$. Under some specific conditions, $c$ may become a
increasing function to $\gamma$ (see Remark~\ref{1.4}
and the graph 1 in each panel in Figure \ref{fig:ct_3_12_allm}).
However, it is not clear if the function $c$ has monotonicity generically. Using
the Newton's method, we solve the system of
equations (\ref{1.19.1}) and (\ref{1.19.2}) and obtain the graph
of $c$ and $t$, respectively. We first consider three ion species
with coefficients satisfying $b_1=1$, $b_2=2$ and $b_1\beta_1 +
b_2\beta_2=a_1\alpha_1=1.2$. Specific values of $\beta_1$ and
$\beta_2$ can be chosen as follows:
\begin{enumerate}
\item[I.]~~~$(\beta_1,\beta_2) =(1.199,0.0005)$,
\item[II.]~~$(\beta_1,\beta_2) =(0.002,0.599)$.
 \end{enumerate}
For each $(\beta_1,\beta_2)$, graphs of $c$ and $t$ corresponding
to the cases of $(a_1, \alpha_1)=(1,1.2)$, $(2,0.6)$ and $(3,0.4)$
are plotted in Figure~\ref{fig:ct_3_12_allm}, respectively.
As for
Theorem~\ref{2-thm5}~(ii)-(iii), our numerical results indicate
that $|c(\gamma)| < t(\gamma)$ for all $\gamma>0$; both
$c(\gamma)$ and $t(\gamma)$ tend to zero as $\gamma $ goes to
infinity. For each fixed $\gamma>0$, the value $t(\gamma)$
increases but the value $c(\gamma)$ decreases as $\alpha_1$
increases. Similar results can also be observed for four ion species with
coefficients satisfying the following conditions:
\begin{enumerate}
\item[]\hspace{-2.0em}Case~1.~~$a_1\alpha_1=\beta_1+2\beta_2+3\beta_3=1.5$,
$(\beta_1,\beta_2,\beta_3)=(0.25,0.25,0.25)$, \\
$a_1=1,2,3,4$, i.e., $\alpha_1=1.5,0.75,0.5,0.375$,
\item[]\hspace{-2.0em}Case~2.~~$\alpha_1+2\alpha_2=\beta_1+2\beta_2=1.5$,
$(\beta_1,\beta_2)=(0.75,0.375)$, \\
$(\alpha_1,\alpha_2)=(0.3,0.6), (0.5,0.5), (0.75,0.375)$,
\item[]\hspace{-2.0em}Case~3.~~$\alpha_1+2\alpha_2=\beta_1+2\beta_2=1.5$,
$(\beta_1,\beta_2)=(0.5,0.5)$, \\
$(\alpha_1,\alpha_2)=(0.3,0.6), (0.5,0.5), (0.75,0.375)$,
\item[]\hspace{-2.0em}Case~4.~~$\alpha_1+2\alpha_2=\beta_1+2\beta_2=1.5$,
$(\beta_1,\beta_2)=(0.3,0.6)$, \\
$(\alpha_1,\alpha_2)=(0.3,0.6), (0.5,0.5), (0.75,0.375)$.
\end{enumerate}
The profiles of $c$ and $t$ associated with Case~1--4 are sketched
in Figure~\ref{fig:ct_4_all}, I--IV, respectively. As for
Figure~\ref{fig:ct_3_12_allm}, various $\alpha_j$'s may result in
different profiles of function $c=c(\gamma)$. However, until now,
all our results only show that the function $c$ is of monotone
increasing or decreasing. This motivates us to see if the
function $c$ becomes a non-monotone function under the other
conditions of $\alpha_i$'s and $\beta_j$'s.

\begin{figure}[hbt]
 \centering{%
  \begin{tabular}{@{\hspace{-0pc}}c}
  \includegraphics[width=3.8in]{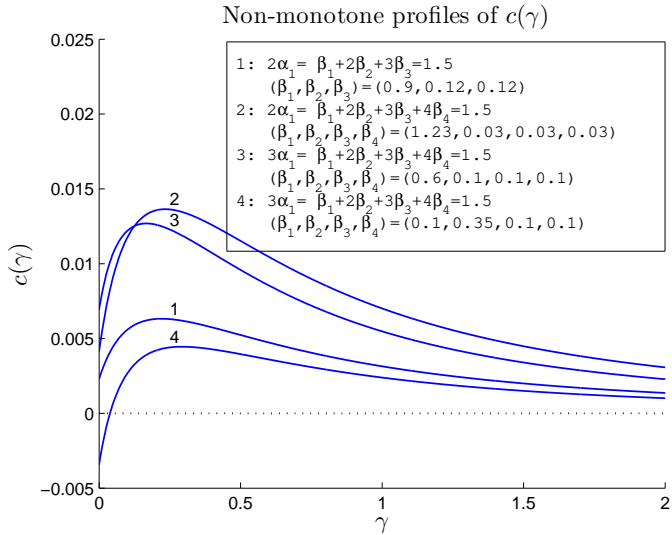}
  \end{tabular}}
 \caption{\small Non-monotone profiles of $c(\gamma)$.
 }\label{fig:ct_5_14_all_z}
\end{figure}

As shown in both Figure~\ref{fig:ct_3_12_allm} and~\ref{fig:ct_4_all}, we observe
that $c(\gamma)$ converges to zero as $\gamma$ goes to infinity.
This is consistent with the results of Theorem~\ref{2-thm5}.
Moreover, the profile of function $c$ can be changed from monotone
decreasing to increasing. Such a behavior of $c$ and the nonlinearity of
equations~(\ref{1.19.1}) and (\ref{1.19.2}) let us believe that
the non-monotone profile of function $c$ may exist. To get the
non-monotone profile of function $c$, we consider the
following conditions:
\begin{enumerate}
\item[A.]~~$2\alpha_1=\beta_1+2\beta_2+3\beta_3=1.5$,
$(\beta_1,\beta_2,\beta_3)=(0.9,0.12,0.12)$,
\item[B.]~~$2\alpha_1=\beta_1+2\beta_2+3\beta_3+4\beta_4=1.5$,
$(\beta_1,\beta_2,\beta_3,\beta_4)=(1.23,0.03,0.03,0.03)$,
\item[C.]~~$3\alpha_1=\beta_1+2\beta_2+3\beta_3+4\beta_4=1.5$,
$(\beta_1,\beta_2,\beta_3,\beta_4)=(0.6,0.1,0.1,0.1)$,
\item[D.]~~$3\alpha_1=\beta_1+2\beta_2+3\beta_3+4\beta_4=1.5$,
$(\beta_1,\beta_2,\beta_3,\beta_4)=(0.1,0.35,0.1,0.1)$.
\end{enumerate}
The non-monotonic profiles of function $c$ with respect to
conditions~A-D are provided in Figure~\ref{fig:ct_5_14_all_z},
1-4, respectively. However, the profiles of functions $t$ and $t-c$ are still monotonically decreasing.

\section{Conclusion}\label{2-mpnpsec:6}

For the binary mixture of
monovalent anions and cations, although CCPB and PB can have very different
solutions with different boundary conditions and other constraints, the
solutions of CCPB equations have
very similar asymptotic behavior as those of PB equations when the
global electroneutrality (\ref{g-neutral}) holds (cf.~\cite{LHLL}).

Situation becomes more complicated in the presence of mixtures of multiple (more than three) species with multivalences. In this paper, we again consider the situations under
global electroneutrality, but the general mixture of
multi-species ions. The (more rigorous) CCPB shows very different asymptotic behaviors
to PB equations under Robin type boundary conditions with various coefficients $\eta_\epsilon$'s.

In particular, the solution $\phie$ of CCPB equation may tend to a
constant $c$ at interior points, and $\pm t$ at boundary points as
$\epsilon$ goes to zero. As $\eta_\epsilon\sim\gamma\epsilon$,
both $t$ and $t-c$ are monotone decreasing functions of $\gamma$.
Physically, $\gamma$ can be regarded as the ratio of the
Stern-layer width to the Debye screening length. Various
conditions can be found theoretically and numerically such that
the function $c$ of $\gamma$ becomes monotone decreasing,
increasing and non-monotonic. While for PB
equation, the solution $\phie$ only tend to zero at interior
points which is independent to $\gamma$. This constitutes one of the
main differences of PB and CCPB equations.

This work is one of our first attempts in systematically studying the ionic fluids.
Much works are needed in the future. In particular, the theoretical justification of the
interesting behavior of $c(\gamma)$ with respect to $\gamma$ under different physical
conditions. The problems involving  multiple spatial dimension domains are for certain to
provide more interesting phenomena of the solutions and also more technical challenges.
Overall, our results again demonstrate that the CCPB equation being a more physical and suitable
model for future applications involving the mixture of multi-species ions.


\section*{Appendix}

For the convenience of the readers, we will list out our previous results for 2 mono-valence species with charges of
opposite signs  situations ~\cite{LHLL} .

Considering CCPB equation~(\ref{2-eqn1}) with $N_1=N_2=1$, $a_1=b_1=1$, in~\cite{LHLL}, we had
established the following results:

\noindent\textbf{(a)} In the electroneutral case ($\alpha_1=\beta_1$):
\begin{enumerate}
\item[(a1)] If
$\lim_{\epsilon\downarrow0}\frac{\epsilon}{\eta_\epsilon}=0$, the
solution $\phie$ approaches zero in $[-1,1]$ as
$\epsilon\downarrow0$. However, $\phie$ has slope of order
$O(1/\eta_\epsilon)$ on the boundary.
\item[(a2)] When $\frac{\epsilon}{\eta_\epsilon}\geq\,C$ for some
positive constant $C$ independent of $\epsilon$, the solution
$\phie$ possesses boundary layers with thickness $\epsilon$.
\end{enumerate}

\noindent\textbf{(b)} In the non-electroneutral case ($\alpha_1\neq\beta_1$):\\
The solution $\phie$ has boundary layers with thickness
$\epsilon^2$ and $\phie(x)-\phie(\pm1)$ tends to infinity with the
leading order term $\log(\epsilon^{-2})$ as $\epsilon\downarrow0$
for $ x\in (-1,1)$. The
values $\phie(\pm 1)$ can be estimated as follows:
\begin{enumerate}
\item[(b1)] If $\frac{\eta_\epsilon}{\epsilon^2}\leq C$, $\phie(1)$
and $\phie(-1)$ converge to different finite values as
$\epsilon\downarrow0$, where $C$ is a positive constant
independent of $\epsilon$.
\item[(b2)] If
$\lim_{\epsilon\downarrow0}\frac{\eta_\epsilon}{\epsilon^2}=\infty$,
both $\phie(1)$ and $\phie(-1)$ diverge to $\infty$, but
$|\phie(1)-\phie(-1)|$ converges to zero as $\epsilon\downarrow0$.
\end{enumerate}
\noindent\textbf{(c)} The difference between the solutions to the CCPB equation~(\ref{2-eqn1}) and the PB equation~(\ref{2-eqnw}) can
be stated as follows:
\begin{enumerate}
\item[(c1)] When $\alpha_1=\beta_1$, the solution of the CCPB equation
(\ref{2-eqn1}) may converge to the solution of the PB equation
(\ref{2-eqnw}). Namely, in the case of
$\alpha_1=\beta_1$, the solution of the CCPB equation (\ref{2-eqn1})
has the same asymptotic behavior as that of the PB equation
(\ref{2-eqnw}).
\item[(c2)] When $\alpha_1\neq\beta_1$, the solution of
the PB equation (\ref{2-eqnw}) remain bounded for $\epsilon>0$. However, as $\alpha_1\neq\beta_1$, the solution
of the CCPB equation (\ref{2-eqn1}) may tend to infinity as
$\epsilon$ goes to zero (see~\textbf{(b)}). This may provide the
difference between the solutions to the CCPB
equation~ (\ref{2-eqn1}) and the PB equation~(\ref{2-eqnw}).
\end{enumerate}



\begin{thebibliography}{}

          %


          %
\bibitem{A1} D. Andelman, {\em Electrostatic Properties of Membranes: The Poisson Boltzmann Theory, Handbook of Biological Physics}, Volume~1, edited by R. Lipowsky and E. Sackmann, Elsevier Science B.V., 603-641, 1995.

\bibitem{APP} F. Andrietti, A. Peres and R. Pezzotta, {\em Exact solution of the unidimensional Poisson-Boltzmann equation for a 1:2 (2:1) electrolyte}, Biophys J. 9:1121-4, 1976.

\bibitem{BNVHEG07} D. Boda, W. Nonner, M. Valisko , D. Henderson, B. Eisenberg, and D. Gillespie, {\em Steric Selectivity in Na Channels Arising from Protein Polarization and Mobile Side Chains}, Biophys. J., 93, 1960--1980, 2007.

\bibitem{BVHEG07} D. Boda, M. Valisko , D. Henderson, B. Eisenberg, and D. Gillespie, {\em Ionic selectivity in L-type calcium channels by electrostatics and hard-core repulsion}, J. Gen. Physiol., 133(5), 497--509, 2009.

\bibitem{BrBr93} C. M. A. Brett and A. A. O, Brett, {\em Electrochemistry. Principles, Methods, and Applications},
Oxford Science Publications, Oxford, 1993.

\bibitem{BCB} M.Z. Bazant, K.T. Chu and B. J. Bayly, {\em Current-Voltage relations for electrochemical thin films}, SIAM J. Appl. Math., 65(5), 1463--1484, 2005.

\bibitem{BCE} V. Barcilon, D. P. Chen, R. S. Eisenberg, and J. W. Jerome, {\em Qualitative properties of steady-state Poisson-Nernst-Planck systems: perturbation and simulation study}, SIAM J. Appl. Math., 57(3), 631--648, 1997.

\bibitem{BCE1} V. Barcilon, D. P. Chen, and R. S. Eisenberg, {\em Ion flow through narrow membrane channels: part II}, SIAM J. Appl. Math., 52(5), 1405--1425, 1992.

\bibitem{BaFa01} A. J. Bard and L. R. Faulkner, {\em Electrochemical Methods}, John Wiley \& Sons, Inc., New York, NY, 2001.


\bibitem{CBW-10} Z. Chen, N.A. Baker and G.W. Wei, {\em Differential geometry based solvation model I: Eulerian formulation}, Journal of Computational Physics, 229, 8231--8258, 2010.

\bibitem{CE1}  D. Chen, P. Kienker, J. Lear and B. Eisenberg, {\em PNP Theory fits current-voltage (IV) relations of a synthetic channel in 7 solutions}, Biophys. J., 68:A370, 1995.

\bibitem{Ch05}  K. T. Chu, {\em Asymptotic Analysis of Extreme Electrochemical Transport}, PhD thesis, NMIT, 2005.

\bibitem{DC} S. Das and S. Chakraborty, {\em Effect of Conductivity Variations within the Electric Double Layer on the Streaming Potential Estimation in Narrow Fluidic Confinements}, Langmuir, 26(13), 11589--11596 2010.

\bibitem{FT} A. Friedman and K. Tintarev, {\em Boundary asymptotics for solutions of the Poisson-Boltzmann equation}, J. Differential Equations, 69 , 15--38, 1987.

\bibitem{GP} A. J. M. Garrett and L. Poladian, {\em Refined derivation, exact solutions and singular limits of
the Poisson Boltzmann equation}, Ann. Phys., 188, 386--435, 1988.

\bibitem{GYW07} W. Geng, S. Yu, and G. Wei, {\em Treatment of charge singularities in implicit solvent models}, J. Chem. Phys., 127, 114106, 2007.

\bibitem{G1} S. Glasstone, {\em An introduction to Electrochemistry}, D. Van Nostrand Company, Inc., Princeton, N.J., 1942.

\bibitem{GW}  H. T. Gordon and J. H. Welsh, {\em The role of ions in axon surface reactions to toxic organic compounds}, J. Cell. Comp. Physiol, 31, 395-419, 1948.

\bibitem{H1} B. Hille, {\em Ion channels of excitable membranes, 3rd
Edition, Sinauer Associates}, Inc., 200).

\bibitem{H2}  V. R. T. Hsu, {\em Almost Newton method for large flux steady-state of 1D Poisson-Nernst-Planck equations}, J. Comput. Appl. Math., 183(1), 1--15, 2005.

\bibitem{HL} J. P. Hsu and B. T. Liu, {\em Current Efficiency of Ion-Selective Membranes: Effects of Local Electroneutrality and Donnan Equilibrium}, J. Phys. Chem. B, 101, 7928--7932, 1997.

\bibitem{Hunter} R. J. Hunter, {\em Zeta Potential in Colloid Science}, Academic Press Inc., 198).

\bibitem{Hunter1} R. J. Hunter, {\em Foundations of Colloid Science}, Oxford University Press, London, 2001.

\bibitem{HyonMori} Y. Hyon, {\em A Mathematical Model For Electrical Activity in Cell Membrane: Energetic Variational Approach}, in preparation, 2015.

\bibitem{I} J. Israelachvili, {\em Intermolecular and Surface Forces}, Academic Press, London, 1992.

\bibitem{K0} J. Keener and J. Sneyd, {\em Mathematical Physiology}, Springer-Verlag, New York, Inc, 1998.

\bibitem{KLL-11} V.Y. Kiselev, M. Leda, A.I. Lobanov, D. Marenduzzo,
and A.B. Goryachev, {\em Lateral dynamics of charged lipids and peripheral proteins in spatially heterogeneous membranes: Comparison of continuous and Monte Carlo approaches}, J. Chem. Phys. 135, 155103, 2011.

\bibitem{K1} G. Kort\"{u}m, {\em In {\it Treatise on Electrochemistry}, 389--394}, London/New York: Elsevier, 1965.

\bibitem{LMB} D. Lacoste, G.I. Menon, M.Z. Bazant, and J.F. Joanny, {\em Electrostatic and electrokinetic contributions to the elastic moduli of a driven membrane}, Eur. Phys. J. E, 28, 243--264, 2009.

\bibitem{LHLL} C. C. Lee, H. Lee, Y. Hyon, T. C. Lin and C. Liu, {\em New Poisson-Boltzmann Type Equations: One-Dimensional Solutions}, Nonlinearity, 24, 431--458, 2011.

\bibitem{Lcc} C. C. Lee, {\em The Charge Conserving Poisson-Boltzmann Equations: Existence, Uniqueness and Maximum Principle}, J. Math. Phys., 55, 051503, 2014.

\bibitem{LC} M. Lee and K. Y. Chan, {\em Non-neutrality in a charged slit pore}, Chem. Phys. Letts., 275, 56--62, 1997.

\bibitem{Lw} W. Liu, {\em One-dimensional steady-state Poisson-VNernst-VPlanck
systems for ion channels with multiple ion species}, J. Differential Equations, 246 428--451, 2009.

\bibitem{Lw2} {\sc W. Liu}: {\em Geometric singular perturbation approach to steady-state Poisson-Nernst-Planck systems}, SIAM J. Appl. Math. 65 (2005), no. 3, 754-766.

\bibitem{MC} J. E. Marsden and A. J. Chorin, {\em A Mathematical Introduction To Fluid Mechanics}, Springer, 1993.

\bibitem{Mori1} Y. Mori, J.W. Jerome, and C.S. Peskin, {\em A Three-dimensional Model of Cellular Electrical Activity}, Bulletin of the Institute of Mathematics Academia Sinica, 2(2), 367--390, 2007.

\bibitem{Mori2} Y. Mori, {\em A Three-Dimensional Model of Cellular Electrical Activity}, PhD thesis, New York University, 2007.

\bibitem{MPR}  W.E. Morf, E. Pretsch and N.F. de Rooij, {\em Theoretical treatment and numerical simulation of potential and concentration profiles in extremely thin non-electroneutral membranes used for ion-selective electrodes}, Journal of Electroanalytical Chemistry 641, 45--56, 2010.

\bibitem{Ne91} J. Newman, {\em Electrochemical Systems}, Prentice-Hall Inc., Englewood Cliffs, NJ, Second edition, 1991.

\bibitem{NP}  B. W. Ninham and V. A. Parsegian, {\em Electrostatic potential between surfaces bearing ionizable groups in ionic equilibrium with physiologic saline solution}, J. Theor. Biol., 31(3), 405--428, 1971.

\bibitem{O} L.H. Olesen, {\em AC Electrokinetic micropumps}, thesis, Department of Micro and Nanotechnology Technical University of Denmark, 2006.

\bibitem{PJ} J. H. Park and J. W. Jerome, {\em Qualitative properties of steady-state Poisson-Nernst-Planck systems: mathematical study}, SIAM J. Appl. Math., 57(3), 609--630, 1997.

\bibitem{RWL} E. Riccardi, J. C. Wang  and A. I. Liapis, {\em Porous Polymer Absorbent Media Constructed by Molecular Dynamics Modeling and Simulations: The Immobilization of Charged Ligands and Their Effect on Pore Structure and Local Nonelectroneutrality}, J. Phys. Chem. B, 113, 2317--2327, 2009.

\bibitem{RABI-04} B. Roux, T. Allen, S. Berneche and W. Im, {\em Theoretical and computational models of biological ion channels}, Quarterly Reviews of Biophysics, 37(1), 15--103, 2004.

\bibitem{R} I. Rubinstein, {\em Counterion condensation as an exact limiting property of solutions of the Poisson-Boltzmann equation}, SIAM J. Appl. Math. 46, 1024--1038, 1986.

\bibitem{R1} R. Ryham, {\em An Energetic Variational Approach To Mathematical Modeling Of Charged Fluids: Charge Phases, Simulation And Well Posedness}, thesis, Pennsylvania State University, 2006.

\bibitem{RLZ} R. Ryham, C. Liu, and L. Zikatanov, {\em Mathematical Models for the Deformation of Electrolyte Droplets}, Discrete Contin. Dyn. Syst. Ser. B, 8(3), 649--661, 2007.

\bibitem{RLW} R. Ryham, C. Liu and Z.Q. Wang, {\em On electro-kinetic fluids: one dimensional configurations}, Discrete Contin. Dyn. Syst. Ser. B, 6(2), 357--371, 2006.

\bibitem{SGN} A. Singer, D. Gillespie, J. Norbury and R.S. Eisenberg, {\em Singular perturbation analysis of the steady-state Poisson-Nernst-Planck system: applications to ion channels}, European J. Appl. Math., 19(5), 541--560, 2008.

\bibitem{WXL13} L. Wan, S. Xu, M. Liao, C. Liu and P. Sheng, {\em Self-consistent approach to global charge neutrality in electrokinetics: A surface potential trap model}, Phys. Rev. X 4, 011042, 2014.

\bibitem{ZMB} F. Ziebert, M. Z. Bazant and David Lacoste, {\em Effective zero-thickness model for a conductive membrane driven by an electric field}, Phys Rev E 81, 031912(1-13), 2010.

\bibitem{ZW-11} Q. Zheng and G.W. Wei, {\em Poisson-Boltzmann-Nernst-Planck model}, J. Chem. Phys. 134, 194101, 2011.

\bibitem{ZWL} S. Zhou, Z. Wang and Bo Li, {\em Mean-Field description of ionic size effects with non-Uniform ionic sizes: A numerical approach}, Phys. Rev. E 84, 021901, 2011.

          \end{thebibliography}
          \end{document}